\documentclass[DIV=13,10pt,a4paper,abstract=true]{scrartcl} 
\addtokomafont{disposition}{\rmfamily}

\usepackage[utf8]{inputenc}
\usepackage[british]{babel}
\usepackage[english=british]{csquotes}
\usepackage[section]{placeins} 
\usepackage{lmodern}
\usepackage{hyperref}
\usepackage{hyperxmp}
\usepackage{xcolor}
\usepackage{tikz}
\usepackage{pgfplots}\pgfplotsset{compat=1.16}\usepgfplotslibrary{patchplots}
\usepackage{wrapfig} 
\usepackage{float} 
\usepackage{booktabs}
\usepackage{mathtools} 
\usepackage{amsmath}
\usepackage{amssymb}
\usepackage{amsfonts}
\usepackage[capitalize,nameinlink,noabbrev]{cleveref}
\usepackage[font={small,it}]{caption}
\usepackage{subcaption}
\usepackage[amsmath,amsthm,thmmarks,hyperref]{ntheorem}
\usepackage
[
	backend=bibtex,
	language=british,
	sorting=none,
    sortcites=true,
	sortlocale=en-GB,
    giveninits=true,
    maxnames=4,
    isbn=true,
    doi=true,
    autocite=superscript
]{biblatex}
\DeclareCiteCommand{\supercite}[\mkbibsuperscript]{%
\iffieldundef{prenote}{}{\BibliographyWarning{Ignoring prenote argument}}%
\iffieldundef{postnote}{}{}}%
{\bibopenbracket\usebibmacro{citeindex}\usebibmacro{cite}\usebibmacro{postnote}\bibclosebracket}{\supercitedelim}{}

\colorlet{internallinkcolour}{green!50!black}
\colorlet{externallinkcolour}{red!50!black}
\hypersetup
{
  pdfauthor          = {R. Paul Wilhelm and Matthias Kirchhart},
  pdftitle           = {An Interpolating Particle Method for the Vlasov–Poisson Equation},
  pdfencoding        = {unicode},
  pdflang            = {en-GB},
  pdfdisplaydoctitle = {true},
  pdfkeywords        = {Vlasov–Poisson equation; Lagrangian method; particle method, 
					   Reproducing kernel Hilbert spaces, radial basis functions},
  unicode            = {true},
  colorlinks         = {true}, 
  allcolors          = {internallinkcolour},
  urlcolor           = {externallinkcolour},
  frenchlinks        = {false},
  pdfborder          = {0 0 0},
  naturalnames       = {false},
  hypertexnames      = {false},
  breaklinks         = {true}
}

\theorempreskip{1ex plus .25ex minus .1ex}
\theorempostskip{1ex plus .25ex minus .1ex}

\theoremstyle{nonumberplain}
\theoremheaderfont{\normalfont\scshape}
\theorembodyfont{\normalfont\itshape}
\theoremseparator{.}
\theoremsymbol{}

\theoremheaderfont{\normalfont\itshape}
\theorembodyfont{\normalfont\small}

\newcommand{\R}{\ensuremath{\mathbb{R}}} 
\newcommand{\N}{\ensuremath{\mathbb{N}}} 

\newcommand{\dx}{\ensuremath{{\mathrm d}x}} 
\newcommand{\dy}{\ensuremath{{\mathrm d}y}} 
\newcommand{\dv}{\ensuremath{{\mathrm d}v}} 
\newcommand{\dt}{\ensuremath{{\mathrm d}t}} 

\newcommand{\ddt}{\ensuremath{\frac{{\mathrm d}}{{\mathrm d}t}}}

\newcommand{\reals}{\ensuremath{\mathbb{R}}}
\newcommand{\integers}{\ensuremath{\mathbb{Z}}}
\newcommand{\naturals}{\ensuremath{\mathbb{N}}}
\newcommand{\spdim}{\ensuremath{\mathrm{d}}} 

\newcommand{\vmax}{\ensuremath{v_{\mathrm{max}}}}

\usepackage{microtype}
\begin{document}

\title{An Interpolating Particle Method for the Vlasov--Poisson Equation}
\author{R.~Paul~Wilhelm\thanks{Applied and Computational Mathematics,
RWTH Aachen University, Schinkelstra{\ss}e~2, 52062~Aachen, Germany. E-mail:~%
\href{mailto:wilhelm@acom.rwth-aachen.de}{wilhelm@acom.rwth-aachen.de}}
\and
Matthias Kirchhart\thanks{E-Mail:~%
\href{mailto:kirchhart@acom.rwth-aachen.de}{kirchhart@acom.rwth-aachen.de}.
The second author has been funded by the German Research Foundation (DFG),
project number 432219818, \enquote{Vortex Methods for Incompressible Flows}.}}
\date{\vspace{-1cm}}
\maketitle

\begin{abstract}
In this paper we present a novel particle method for the Vlasov--Poisson
equation. Unlike in conventional particle methods, the particles are not
interpreted as point charges, but as point values of the distribution
function. In between the particles, the distribution function is reconstructed
using mesh-free interpolation. Our numerical experiments confirm that
this approach results in significantly increased accuracy and noise reduction.
At the same time, many benefits of the conventional schemes are preserved.
\end{abstract}

\section{Introduction}
The Vlasov--Poisson system is a simplified model for the evolution of plasmas in
their collisionless limit, as they occur in, for example, nuclear fusion
devices. In dimensionless form this system is given by:
\begin{gather}
\label{eqn:vlasov}
\partial_t f + v\cdot\partial_x f - E\cdot\partial_v f = 0, \\
\label{eqn:Edef}
E\coloneqq -\nabla\varphi,\\
\label{eqn:poisson}
-\Delta\varphi = \rho,\\
\label{eqn:density}
\rho(t,x)\coloneqq 1 - \int_{\reals} f(t,x,v)\,\dv.
\end{gather}
Here, $f=f(t,x,v)$ is the electron distribution function, i.\,e., $f(t,x,v)$
describes the probability density of electrons having velocity $v\in\reals$
and location $x\in\reals$ at time $t\in\reals$.
We will assume that $f$ is periodic in $x$ with period $L>0$, i.\,e., $f(t,x,v)=f(t,x+kL,v)$ 
for any $x\in\reals$ and $k\in\integers$. Therefore it suffices to look at $x\in [0,L]$.
We need to demand that $f$ is normalised 
such that:
\begin{equation}\label{eqn:fnormalization}
\frac{1}{L}\int_{[0,L]}\int_{\reals} f(t,x,v)\,\dv\,\dx
= 1.
\end{equation}
\Cref{eqn:density} defines the charge density $\rho$, were the additional
\enquote{1} stems from the assumption of a uniform ion-background and thus 
ensures overall neutrality 
\begin{equation}
\label{eqn:neutrality_rho}
\int_0^L\rho(t,x)\,\dx = 0. 
\end{equation}
Neglecting collisions and the magnetic field, the Vlasov equation~\eqref{eqn:vlasov} 
then describes the evolution of $f$ under the influence of the self-consistent  
electrical field $E=E(t,x)$, given in terms of the electric potential 
$\varphi=\varphi(t,x)$, which in turn is  given as the solution of the Poisson
equation~\eqref{eqn:poisson}.

Particle-in-Cell methods (PIC) are a long established tool to obtain numerical
approximations to solutions of this system. However, it is well-known that
these methods suffer from \enquote{numerical noise} and have a low convergence
order\autocite{raviart_analysis_1985,ameres_stochastic_2018}. 
For this reason there has been an increased interest in high-order
grid-based methods\autocite{filbet_comparison_2001, filbet_conservative_2001}. 
These methods often have good stability
properties and generalise to arbitrary order. Unlike particle methods, however,
they introduce numerical dissipation; especially when the true solution develops
features that are smaller than the grid size. This is also the case for 
so-called \emph{remapped} or \emph{remeshed} particle methods:
here the particles are periodically remapped onto a Cartesian grid to avoid 
the aforementioned \enquote{numerical noise}%
\autocite{wang_particle--cell_2011}. However, this remeshing effectively acts
as a low-pass filter, smearing out features below the grid resolution.

Particle methods without remapping, on the other hand, are based on an analytic
solution: particles simply follow the characteristic lines, and are thus free
of numerical dissipation. It has recently been shown that the \enquote{numerical
noise} is actually the result of interpreting the particle field as a quadrature
rule\autocite{kirchhart_particle_2019}. Instead of interpreting a particle
field as a set of points with associated \emph{weights}, it should be
interpreted as a set of points with associated \emph{function values}. Instead
of \emph{regularising} a quadrature rule, one should try to \emph{interpolate}
between the known function values. In this work we want to show that this in
fact leads to particle methods that achieve accuracies comparable to grid-based
methods, without the associated introduction of spurious numerical dissipation.
Nonetheless, as our numerical experiments will show, aliasing effects limit the 
accuracy of these approaches, giving rise to phenomena which have not yet been 
reported for particle methods. 

\section{Related Literature and Methods}
The Vlasov--Poisson--Maxwell equations were introduced by Vlasov in his 1938
seminal paper\autocite{vlasov_vibrational_1938}. The Vlasov--Poisson system
\eqref{eqn:vlasov}--\eqref{eqn:density} results when magnetic effects are
neglected. In 1945 Landau gave a first analysis of a linearisation of this
system close to an equilibrium state\autocite{ter_haar_61_1965}. 
Arsen'ev gave the first regularity and well-posedness analysis for the 
one-dimensional Vlasov--Poisson system\autocite{ARSENEV1975131,1975ZVMMF}. Ukai
and  Okabe extended the results to the two-dimensional 
case\autocite{ukai_classical_1978}. Pfaffelmoser gave an existence and
uniqueness result for classical initial data in the three-dimensional case.
Lions and Perthame gave an result for weak
data~\autocite{pfaffelmoser_global_1992,lions_propagation_1991}.

Mouhot and Villani extended the analysis done by Landau to the non-linear
system  in their seminal paper\autocite{mouhot_landau_2011}. We refer to their
work for more details; it also contains an extensive review of theory,
bibliography and historic remarks.  

A general overview of numerical methods for the Vlasov--Poisson equation is
provided in the books by Birdsall, Langdon and Glassey%
\autocite{birdsall_plasma_1985,glassey_cauchy_1996}. Filbet and
Son\-nen\-dr\"{u}cker compared different Eulerian
approaches\autocite{filbet_comparison_2001,filbet_conservative_2001}.

Particle methods have their origins in fluid dynamics%
\autocite{rosenhead1931,Harlow1,HarlowEvans2}. An introduction to particle
methods and their application to the Vlasov--Poisson equation can be found in
the articles by Raviart and Cottet%
\autocite{raviart_analysis_1985,cottet_particle_1984} as well as in Hockney's
book\autocite{hockney_computer_2021}.

\emph{Particle-in-Cell} methods (PIC) were originally developed by Evans and 
Harlow for applications in hydrodynamics\autocite{HarlowEvans2}. A literature 
overview can be found in Hockney and Eastwood's book\autocite{hockney_computer_2021}. 
Denavit suggested to use remeshing to reduce particle noise\autocite{DENAVIT197275}
and Wang, Miller, Colella, Myers and Straalen built conservative high order PIC methods 
with remeshing\autocite{wang_particle--cell_2011,myers_4th-order_2016}. 
Cottet and Raviart gave an analysis of the PIC method for the Vlasov--Poisson 
equation\autocite{pic_cottetRaviart}.
A more statistical approach to PIC and a combination with Monte-Carlo based methods is 
presented by Ameres in his recent PhD thesis\autocite{ameres_stochastic_2018}. Ameres
also discusses particle  noise and its influence on the convergence of
(statistical) PIC methods. 

Semi-Lagrangian schemes were proposed by Rossmanith and Seal, Sonnendrücker and
Besse as well  as Charles, Després and Mehrenberger%
\autocite{besse_semi-lagrangian_2003,charles_enhanced_2013,ROSSMANITH20116203}. 
An extension of semi-Lagrangian schemes to higher dimensions and comparisons to 
other approaches were presented by Cottet\autocite{cottet:hal-01584107}. 

Related to the approach we propose, Russo and Strain suggested using 
interpolation for a purely Lagrangian scheme in the context of vortex 
methods\autocite{RUSSO1994291}.
However, unlike our scheme, their method requires the generation of a 
triangular mesh using the particle locations as mesh nodes in every time step.
Triangulation is expensive and in particular the method does not scale well
with higher dimensions. In contrast, our apprach is based on ideas from 
mesh-free methods.

A general overview of the \emph{reproducing kernel Hilbert space} framework
(RKHS) is given in the books by Wendland and Fasshauer%
\autocite{wendland_scattered_2004,fasshauer_meshfree_2007}. For brevity we will
refer to methods using the RKHS framework as kernel-based methods.  An analysis
of the stability of kernel-based interpolation in Sobolev spaces can be found
in the article by de Marchi and Schaback%
\autocite{de_marchi_stability_2010} as well as Rieger's PhD-thesis%
\autocite{rieger_sampling_2009}. The books of Fasshauer and Wendland also give
an overview of efficient  implementation of techniques for these methods. 

Reproducing kernels are used in Eulerian-based approaches for transport
equations where this idea was introduced amongst others by Schaback and
Franke\autocite{franke_solving_1998}. A kernel-based interpolation approach in
the Semi-Lagrangian framework was proposed by Iske and Behrends%
\autocite{behrens_grid-free_2002}. This ansatz was further developed by several
authors for both the linear transport equation and some non-linear equations
like the shallow water equation%
\autocite{iske_radial_2003,hunt_mesh-free_2005,bonaventura_kernel-based_2011,
shankar_mesh-free_2018}. 

Finally we want to mention that the RKHS framework was already used in the
context of  SPH methods. Several authors worked in this context on the
so-called RKHS particle method and also applied it to the Vlasov--Poisson
equation\autocite{liu_reproducing_1995}.

\section{Solution-structure of the Vlasov--Poisson Equation}
\label{section:vp-theory}

While in reality the electric field $E(t,x)$ needs to be computed from 
the unknown function $f$ via \eqref{eqn:Edef}--\eqref{eqn:density}, let
us for the moment assume $E(t,x)$ was given for all times $t$.
In this case \eqref{eqn:vlasov} is a linear transport equation and can
be written as
\begin{align*}
\partial_t f(t,z) + \left( a(t,z) \cdot \nabla_z \right) f(t,z) = 0,
\end{align*}
where $z \coloneqq \left( x,v \right) \in \R \times \R$, $t \ge 0$ and
$ a(t,z) := \left( v , -E(t,x) \right)$.
 
This equation can be solved using the method of characteristics. To this end,
for each initial time $s$ and position $z$, we define the trajectory
$\Phi(t;s,z)$ as the solution of the following initial value
problem:
\begin{gather}
\label{eqn:char1}
\ddt \Phi(t;s,z) = a\left(t,\Phi(t;s,z)\right),  \\
\label{eqn:char2} 
\Phi(s;s,z) = z.
\end{gather}
We will also use the notation $\Phi_s^t(z)\coloneqq\Phi(t;s,z)$. With these
definitions it is a classical result that $\Phi_s^t$ is a well-defined
diffeomorphism with inverse $\Phi_t^s=(\Phi_s^t)^{-1}$. Intuitively, 
$\Phi_s^t(z)$ tells us where a \enquote{particle} at location $z$ at time
$s$ was or will be at a another time $t$.

Using the flow-map the solution of \eqref{eqn:vlasov} can be written as
\begin{equation}
\label{eqn:function_values_f}
f(t,\Phi_0^t(z)) = f_0(z).
\end{equation}

Thus, if we track a finite number of \enquote{particles} $ z_1,\dotsc,z_N$,
using \eqref{eqn:char1}--\eqref{eqn:char2}, we know the value of the
solution $f$  at the current particle positions at any time $t\ge 0$, using
\eqref{eqn:function_values_f}. This is the motivation behind particle methods;
they differ in the way how  values of $f$ are obtained in between the particles.

In the one-dimensional case, for classical initial data satisfying a decay
condition,  Raviart and Cottet have proven, using the work of Ukai and Okabe%
\autocite{ukai_classical_1978}, that \eqref{eqn:vlasov}--\eqref{eqn:density}
has a unique solution\autocite[Theorem~1]{cottet_particle_1984}. Thus our
initial assumption is justified in the sense that the electric field is indeed
well-defined.  In a numerical method, the electric field needs to be computed
from the current approximation of $f$. 

\section{Interpolating Particle Methods}
In this section we will first discuss the general structure of interpolating
particle methods for the Vlasov--Poisson equation. To this end, we will first
give a general algorithm, whose individual steps will be explained in detail
in the subsequent subsections.

\subsection{Overview}
The general structure of interpolating particle methods is as follows:
\begin{enumerate}
\item Subdivide the computational domain into a Cartesian grid of widths
$h_x$ and $h_v$. Take a sample $f_i\coloneqq f_0(x_i,v_i)$, $i=1,\ldots,N$ in
each of the grid's cells. These samples may---but do not need to be---taken at
the respective cell centres.
\item Set $t=0$. Enter the time-step loop:
\begin{enumerate}
\item\label{item:loop} Compute an interpolant $f_{h,\sigma}$ on the current set
of particles $\bigl(x_i(t),v_i(t)\bigr)$ and function values $f_i$,
$i=1,\ldots, N$. 
\item Compute the charge density:
\begin{equation*}
\rho_{h,\sigma}(t,x) \coloneqq
1-\int_{\reals}f_{h,\sigma}(t,x,v)\,\dv.
\end{equation*}
\item Solve the Poisson equation for the electric potential:
\begin{equation*}
-\Delta\varphi_{h,\sigma} = \rho_{h,\sigma},
\end{equation*}
and define the approximate electric field as
$E_{h,\sigma}(t,x)\coloneqq -\nabla\varphi_{h,\sigma}(t,x)$.
\item Advance the following system of ODEs one step $\Delta t$ in time, using,
e.\,g., the symplectic Euler method:
\begin{equation*}
\left\lbrace
\begin{aligned}
\frac{{\mathrm d}x_i}{\dt}(t) &= v_i(t), \\
\frac{{\mathrm d}v_i}{\dt}(t) &= -E_{h,\sigma}\bigl(t,x_i(t)\bigr),
\end{aligned}
\qquad i=1,\dotsc,N.
\right.
\end{equation*}
Note: For higher order methods, one needs to repeat steps 2a--c 
for each stage of the Runge--Kutta method to avoid introducing
splitting errors.
\item Set $t\mapsto t + \Delta t$ and go to \cref{item:loop}.
\end{enumerate}
\end{enumerate}

On the one hand, this algorithm closely mirrors conventional blob-based
methods. The key difference lies in \cref{item:loop}. In a conventional blob
method one would chose some blob-function $\zeta_\sigma(x,v)$ of
blob-size $\sigma>0$ and set $f_{h,\sigma} =
h_x^\spdim h_v^\spdim\sum_{i=1}^{N}f_i\zeta_\sigma(x-x_i,v-v_i)$. The resulting
approximation, however, will usually not interpolate the data and contain large
errors. If, on the other hand, an appropriate
interpolation scheme is employed, drastic improvements in accuracy can be
achieved.

\subsection{Construction of Interpolants}
For any given particle field $(x_i,v_i)$ with associated data $f_i$,
$i=1,\dotsc,N$, there are of course infinitely many possible interpolants. This
gives us the freedom to request further conditions. In our case, we demand the
following:
\begin{itemize}
\item\emph{Accuracy.} The interpolant should converge to the true function $f$
at high order, i.\,e., fulfil error bounds of the shape $\mathcal{O}(h^s)$,
where $h\coloneqq\max\lbrace h_x,h_v\rbrace$ is the particle spacing and $s>0$
is the (hopefully high) convergence order.
\item\emph{Stability.} The interpolant should react gracefully to disturbances
in the data $(x_i,v_i)$ and $f_i$, $i=1,\dotsc,N$.
\item\emph{Efficiency.} Construction and evaluation of interpolants need to be
carried out on computers in a fast manner and should require only little extra
storage. In particular the algorithm should be easily parallelisable.
\item\emph{Ease of integration.} Given an interpolant, it must be possible to
compute the charge density $\rho_{h,\sigma}$ both accurately and efficiently.
\end{itemize}

These constraints are fulfilled by piece-wise, tensorised kernel-interpolants,
which we will describe in more detail. 

\subsubsection{Kernel-based Interpolants}
For brevity, we will again use the abbreviation $z=(x,v)\in\reals\times
\reals$ for coordinates in the phase space. Kernel-based interpolants
are functions of the shape $f_{h,\sigma}(z) =\sum_{i=1}^N c_i k(z,z_i)$, where
$\mathsf{c}=(c_1,\dotsc,c_N)^\top\in\reals^N$ is a coefficient vector, and
$k(\cdot,\cdot)$ is a suitable kernel function. The coefficient vector needs
to be chosen such that the interpolation conditions are fulfilled:
\begin{multline}\label{eqn:kernel-system}
\forall i\in\lbrace1,\dotsc,N\rbrace:\quad f_{h,\sigma}(z_i) = f_i \\
\iff
\underbrace{%
\begin{pmatrix}
k(z_1,z_1) & k(z_1,z_2) & \cdots & k(z_1,z_N) \\
k(z_2,z_1) & k(z_2,z_2) & \cdots & k(z_2,z_N) \\
\vdots     & \vdots     & \ddots & \vdots     \\
k(z_N,z_1) & k(z_N,z_2) & \cdots & k(z_N,z_N)
\end{pmatrix}}_{\eqqcolon\mathsf{K}}
\underbrace{%
\begin{pmatrix}
c_1 \\ c_2 \\ \vdots \\ c_N 
\end{pmatrix}}_{\eqqcolon\mathsf{c}}
=
\underbrace{%
\begin{pmatrix}
f_1\\
f_2\\
\vdots\\
f_N
\end{pmatrix}}_{\eqqcolon\mathsf{f}}.
\end{multline}
For a given kernel-function $k$, the interpolation problem thus reduces to
solving the linear system $\mathsf{Kc=f}$, which can be achieved using standard
methods.

Classical choices of kernels are radial kernels, i.\,e.:
\begin{equation}
k(z,\tilde{z}) = b\left(\frac{|z-\tilde{z}|}{\sigma}\right),
\end{equation}
where $b:\reals_+\to\reals$ is called the \emph{radial basis function}, and
$\sigma>0$  is a scaling parameter which needs to be chosen depending on the
problem, but \emph{independent of $h$} to ensure convergence. Typical choices
for $b$ are Gaussians ($b(r)=\exp{(-r^2)}$), or Wendland's functions%
\autocite[Section~9.3]{wendland_scattered_2004}, which are compactly supported 
piece--wise polynomials $b(r) = b^{W}_{d,n}(r)$, see \Cref{table:wendland_kernel}.
In RKHS literature other often used kernels include
inverse multi--quadratics and thin--plate splines. 
An overview can be found in the books of Wendland and 
Fasshauer\autocite{wendland_scattered_2004,fasshauer_meshfree_2007}.
The appropriate choice of kernel depends heavily on the expected solution 
space. These choices in particular guarantee that the kernel matrix
$\mathsf{K}$ always is symmetric positive definite, such that the system~%
\eqref{eqn:kernel-system} always has a unique solution and can be solved using
the Cholesky decomposition.

\begin{table}
\centering
\begin{tabular}{cc}
\toprule
Function & Formula \\
\midrule
$b^W_{1,2}(r)$ & $(1-r)_+^5(8r^2+5r+1)$   \\
$b^W_{2,2}(r)$ & $(1-r)_+^6(35r^2+18r+3)$ \\
$b^W_{4,2}(r)$ & $(1-r)_+^7(48r^2+21r+3)$ \\
\bottomrule
\end{tabular}
\caption{\label{table:wendland_kernel}Examples of Wendland's radial basis functions
$b^W_{d,n}\in C^{2n}$ where $d\in\N$ is the spatial dimension and $n\in\N$ the
function's order.\autocite[Section~9.3]{wendland_scattered_2004} 
(We deviate from Wendland's notation, who uses $\varphi_{d,n}$ instead.) For brevity,
we write $(1-r)_{+} \coloneqq \max\lbrace0,1-r\rbrace$.}
\end{table}

Note that this approach greatly differs from conventional blob-methods.
Superficially, they both use approximations of the shape
$\sum_{i=1}^N c_i k(z,z_i)$. However, in conventional methods the
resulting approximations usually \emph{do not interpolate} the function values
$f_i$, the coefficient vector \emph{is fixed over time}, the blob width
\emph{$\sigma$ depends on $h$}, and \emph{no linear system} needs to be solved.

It can be shown that these interpolants have several beneficial mathematical
properties, such as minimising the so-called native space 
norms\autocite[Chapter 10]{wendland_scattered_2004}. 
For the Wendland kernels the native spaces are isomorphic and norm-equivalent to
Sobolev spaces, where the regularity of the Sobolev space depends on the
order of the kernel. 
The norm-minimising property guarantees both accuracy and 
stability\autocite[Chapters~10--13]{wendland_scattered_2004}. 
While these interpolants are essentially ideal from the perspective of accuracy
and stability, there are practical hurdles to their application in our
context:
\begin{itemize}
\item Evaluation is costly when the number of particles $N$ is large.
To evaluate $f_{h,\sigma}$ at a single location $z$, it is necessary to 
perform a summation over all points $z_i$, $i=1,\dotsc,N$.
\item The kernel-matrices $\mathsf{K}$ are densely populated and tend to be
extremely ill-conditioned. This excludes the use of iterative solvers.
Especially for Gaussians the matrices quickly become singular within machine
precision. This also excludes the use of \enquote{fast algorithms} like
multipole methods to speed up evaluation and to avoid explicitly storing
$\mathsf{K}$.
\item Integration along the $v$-direction is possible, but difficult.
\end{itemize}

The first two points are only problematic when using very fine discretisations
with large numbers of particles $N$. This problem can be alleviated using
piece-wise interpolants, as described in \cref{subsec:piece-wise}. Integration
becomes significantly easier when one uses tensor--product kernels instead.

\subsubsection{Tensorised Kernels}
Using the Euclidean distance $|\cdot|$, a radial basis function $b$ can be
turned into a kernel in arbitrary high spatial dimensions. An alternative
approach is to use tensorised kernels, which result from multiplying low-%
dimensional kernels
\begin{equation}
k(z,\tilde{z}) = k\bigl((x,v),(\tilde x,\tilde v)\bigr) =
b\left(\frac{|x-\tilde x|}{\sigma_x}\right)
b\left(\frac{|v-\tilde v|}{\sigma_v}\right),
\end{equation}
e.\,g., with $b(r) = b^{W}_{1,2}(r)$ from \Cref{table:wendland_kernel}.
Kernels of this form inherit most of the favourable properties of radial
kernels. In particular, when choosing the same radial basis function $b$ as
above, they also always result in symmetric positive definite kernel matrices
$\mathsf{K}$. They give equally high asymptotic convergence orders, albeit
with different error constants. The only potential drawback is that now
for a single evaluation of $k$ multiple evaluations of $b$ are necessary,
which is, however, only relevant if $b$ is expensive to evaluate.

The main benefit in our case is that such interpolants are significantly easier
to integrate along a single coordinate axis. Assume we are given an interpolant
$f_{h,\sigma}(z)=\sum_{i=1}^{N}c_ik(z,z_i)$ with a tensorised kernel $k$.
We then can compute the charge density $\rho_{h,\sigma}$ as follows:
\begin{equation}
\label{eqn:rho_h_sigma}
\rho_{h,\sigma}(x) = 1 - \int_{-\infty}^{\infty}f_{h,\sigma}\,{\mathrm d}v 
= 1 - \sum_{i=1}^{N} c_i b\left(\frac{|x-x_i|}{\sigma_x}\right)
\underbrace{\int_{-\infty}^{\infty}b\left(\frac{|v-v_i|}{\sigma_v}\right)\,{\mathrm d}v}_{\coloneqq \Lambda},
\end{equation}
where the last integral is a constant that only depends on $b$ and $\sigma_v$.
Once this constant has been computed, evaluation of $\rho_{h,\sigma}$ then
reduces to computing the sum and the evaluation of $b$. In the case of
Wendland's functions, $b$ is a piece-wise polynomial and the integral can be
easily evaluated analytically.

\subsubsection{\label{subsec:piece-wise}Piece-wise Interpolants}
In case of the Vlasov--Poisson equation, only integrals along the $v$-direction
of $f$ are taken. In particular, no derivatives or point evaluations of $f$ are
required. It is thus not necessary to construct a globally smooth interpolant
$f_{h,\sigma}$. This justifies the use of piece-wise interpolants: the
computational domain is divided into a disjoint union of axis-aligned boxes,
each of which containing only a small number of particles $N_{\mathrm{box}}$,
where for each box we demand that
$N_{\mathrm{min}}\leq N_{\mathrm{box}}\leq 2N_{\mathrm{min}}$ for a fixed,
user-defined parameter $N_{\mathrm{min}}$. In our experience, chosing
$50\leq N_{\mathrm{min}}\leq 200$ usually suffices for the one-dimensional
case.

Then in each of these boxes a local, kernel-based interpolant is computed. This
way, the size of the kernel-system \eqref{eqn:kernel-system} remains bounded:
instead of solving one large system of dimension $\reals^{N\times N}$, we now
solve many small systems of maximal dimension $\reals^{2N_{\mathrm{min}}\times
2N_{\mathrm{min}}}$. As $N_{\mathrm{min}}$ is a user-defined constant, the cost 
for solving this local system remains constant as well. Below we will describe
a simple subdivision scheme motivated by $kd$-trees that guarantees
$N_{\mathrm{min}}\leq N_{\mathrm{box}}\leq 2N_{\mathrm{min}}$. Thus, solving all
of these local systems separately, one ends up with an overall optimal complexity
of $\mathcal{O}(N)$. 

Our subdivision is based on $kd$-trees using the so-called cyclic splitting 
rule. To this end, let $Z\subset\R^D$ with $D\in\N$ be a point cloud and let
$B=\reals^D$ denote the initial box. Set $n=1$ and fix a minimal number of
points per box $N_{\mathrm{min}}\in\N$. The algorithm then proceeds as follows:
\begin{enumerate}
\item If $\vert Z\vert<2N_{\mathrm{min}}$ stop and return $B$ and $Z$.
\item Split the box $B$ into two: $B=B_1\cup B_2$, where
\begin{equation}
\begin{split}
B_1 &= \lbrace z\in B\,|\,z_n \leq\text{median $n$-coordinate of the set $Z$}\rbrace,\\
B_2 &= B\setminus B_1.
\end{split}
\end{equation}
Similarly, split the point cloud into two: $Z=Z_1\cup Z_2$, where
$Z_1 = Z\cap B_1$ and $Z_2=Z\cap B_2$.
\item If $n < D$ increase $n$ by 1, else set $n = 1$.
\item Recursively apply this procedure on to $B_1$, $Z_1$ and
$B_2$, $Z_2$.
\end{enumerate}
For more details see Wendland's monograph\autocite[Chapter~14.2]{wendland_scattered_2004}.
Other spatial sub-division schemes are certainly possible, but we found this
simple approach to deliver satisfactory results.

The integration of the resulting interpolant $f_{h,\sigma}$ along the $v$-axis
gets only slightly more complicated. Suppose we are given a finite set of
locations $x^1, x^2, x^3, \dotsc, x^{N_\rho}$ at which we want to evaluate
$\rho_{h,\sigma} = 1 - \int f\,{\mathrm d}v$. We use superscripts to
distinguish these points from the particle locations $z_i=(x_i,v_i)$,
$i=1,\dotsc,N$. This evaluation can be achieved using the following
algorithm:
\begin{enumerate}
\item For $i=1,\dotsc,N_\rho$ set $\rho_{h,\sigma}^i\leftarrow 1$.
\item For each box $B$ of the piece-wise interpolant $f_{h,\sigma}$:
\begin{enumerate}
\item Find the evaluation points $x^{i}$ with
$\bigl(\lbrace x^{i}\rbrace\times\reals\bigr)\cap B\neq\emptyset$.
\item For each such point $x^{i}$ set:
\begin{equation}
\label{eqn:eval_rho_pw_tensor}
\rho_{h,\sigma}^i\leftarrow \rho_{h,\sigma}^i -
\int_{v_{\mathrm{min}}(B)}^{v_{\mathrm{max}}(B)} f_{h,\sigma}(x^i,v)%
\,{\mathrm d}v.
\end{equation}
\end{enumerate}
\end{enumerate}
Here $v_{\mathrm{min}}(B)$ and $v_{\mathrm{max}}(B)$ denote the respective
minimum and maximum $v$ co{\-}ordinates of the axis-aligned box $B$. In the spirit 
of equation \eqref{eqn:rho_h_sigma}, the last integral can be evaluated exactly 
when using tensorised kernels and if the radial function $b$ can be integrated
analytically. This is trivially the case for the piece-wise polynomial Wendland
kernels $b = b^{W}_{1,n}$, $n\in\N$.  Thus, in this case, integration can
be carried out efficiently and exactly. We also remark that this algorithm can be
efficiently parallelised.

In the remainder of the paper, when using piece-wise interpolants, we will
refer  to this as the piece-wise method or the PW method; in contrast, when
using global kernel interpolants we will use the term \enquote{direct method}.

\subsection{Computation of the Electric field}
Given the numerical approximation $\rho_{h,\sigma}$, one needs to solve 
the Poisson equation with periodic boundary conditions to obtain 
$E_{h,\sigma}$, i.\,e., one has to solve:
\begin{gather}
\label{eqn:poisson-num}
-\Delta \varphi_{h,\sigma} = \rho_{h,\sigma} \\
\label{eqn:poisson-bound-cond-num}
\varphi_{h,\sigma}(0) = \varphi_{h,\sigma}(L) = 0.
\end{gather}
From this one can compute $E_{h,\sigma} = -\nabla_x \varphi_{h,\sigma}$. In this
work we use a standard Galerkin method and periodic B-Splines on a uniform grid.

Alternatively it is possible to exploit the solution structure of $f_{h,\sigma}$
and $\rho_{h,\sigma}$, when using tensorised Wendland-kernels.  For the
one-dimensional Poisson equation with periodic boundary conditions the Green's
function $G(x,y)$ and it's derivative $K(x,y)=\partial_x G(x,y)$ are explicitly
known,\autocite{cottet_particle_1984} such that we can write
\begin{equation}
E_{h,\sigma}(x) = \int_0^L K(x,y)\rho_{h,\sigma}(y) \dy.
\end{equation}

Now, using \eqref{eqn:rho_h_sigma}, the above equation turns into
\begin{equation}
\label{eqn:numerical-E-alternative}
E_{h,\sigma}(x) = \int_0^L K(x,y) \dy 
- \Lambda \sum_{i=1}^N c_i\int_0^L K(x,y)b\left(\frac{|x-x_i|}{\sigma_x}\right) \dy.
\end{equation}
Both integrals can be evaluated analytically. This would result in a speed-up
and for a given approximation $f_{h,\sigma}$ give the exact solution of
the electric field. However, our benchmarks showed that the interpolation process 
to obtain $f_{h,\sigma}$ takes several times more computation time than solving
for the electric field $E_{h,\sigma}$, even when going to very high resolution in
the numerical Poisson solver. We therefore decided against implementing equation
\eqref{eqn:numerical-E-alternative}.

\subsection{Remarks on Computational Complexity and Feasibility}
The computationally most expensive steps in both the direct and the piece-wise
approaches is the solution of the kernel systems \eqref{eqn:kernel-system}. When
the Cholesky decomposition is used, this results in complexities of:
\begin{itemize}
\item $\mathcal{O}(\tfrac{1}{3}N^3)$ operations for the direct approach, using
the Cholesky decomposition on a single $\R^{N\times N}$ system.
\item $\mathcal{O}(\tfrac{8}{3}N N_{\mathrm{min}}^2)$, using the Cholesky
decomposition on at most $N/N_{\mathrm{min}}$ systems of dimensions less or
equal to $\R^{2N_{\mathrm{min}}\times2N_{\mathrm{min}}}$.
\end{itemize}

One thus immediately sees that the direct approach quickly becomes infeasible,
while the piece-wise approach achieves $\mathcal{O}(N)$ scaling, with the hidden
constant scaling as $N_{\mathrm{\min}}^2$. As mentioned before, in one spatial
dimension it suffices to chose  $50\leq N_{\mathrm{min}}\leq 200$. Therefore,
the constant $N_{\mathrm{min}}^2$ appears to be rather large.

One should keep in mind, however, that the bottleneck of many modern computer
systems typically is not computational power, but memory bandwidth and latency. 
The solution of the linear systems for each box is a dense, local operation:
highly optimised implementations are readily available and can make optimal
use of the processor's arithmetic units. For this reason, as our experiments
will show, the performance difference compared to conventional PIC is not
as dramatic as one might expect on first sight.

Nevertheless, we expect that in higher dimensions $N_{\mathrm{min}}$ will need to
be chosen larger as well, thereby reducing the method's efficiency. This problem
can likely only be alleviated with suitable preconditioners and iterative 
methods -- an ongoing research topic in the RKHS community.\autocite[Chapter~34]%
{fasshauer_meshfree_2007} This method, in its current form, is best suited for
lower dimensional problems.

\section{Elements of a Convergence Analysis}
\label{section:conv-order}
In the following we will sketch a proof for the theoretic convergence 
order of our method. To this end we will restrict ourselves to the linear case, 
i.\,e., we assume that the electric field $E(t,x)$, and therefore the velocity
field  $a(t,\cdot)$ to be known at all times $t$. We will neglect the 
time-integration error when solving \eqref{eqn:char1} and \eqref{eqn:char2}.
In other words, we assume there is no error in the particle locations, such that
all times $t$ the particle cloud carries the values of the exact solution:
$f\bigl(t,z_i(t)\bigr)=f_i$, $i=1,\dotsc,N$.

Furthermore, we will restrict ourselves to radial instead of tensorised 
kernels. We expect similar results to hold for both radial and tensorised
kernels and numerical experiments support this hypothesis as well. However,
their respective native spaces would differ slightly, thus we would need to
give a technical and lengthy derivation of the correct estimates.

We only consider the direct method. The related convergence result for the 
piecewise method can be proven analogously when interpreting the piecewise 
interpolant as a approximation of the global interpolant, thus having locally 
the same convergence order.

We assume we are given initial data $f_0: \R^2 \rightarrow [0,+\infty)$ which
is smooth enough  and periodic in the first component with period $L>0$.
Furthermore let $f_0$ satisfy
\begin{equation}
\label{eqn:modified_decay_condition}
\exists\vmax>0:\ \forall x\in\R, \vert v\vert\ge v_{\max}:\ 
f_0(x,v) = 0
\end{equation}
and
\begin{equation}
\label{eqn:integral-condition}
\frac{1}{L} \int_0^L\int_{-\infty}^{\infty} f_0(x,v) \dv \dx = 1.
\end{equation}

Define $Q \coloneqq [0,L]\times [-\vmax,\vmax]$ and fix a finite set of 
points $Z_0\subset Q$, $Z_0 = \{ z_1,...,z_N\}$ with the respective function
values $\mathsf{f}=(f_1,\dotsc,f_N)=(f_0(z_1),\dotsc,f_0(z_N)\bigr)\in\reals^N$.
Fix a Wendland kernel $k(z,\tilde{z})=b^{W}_{2,n}(|z-\tilde{z}|/\sigma)$
of order $n\in\N$, $n > 1$, and without loss of generality set the scaling parameter
$\sigma = 1$. We define the fill distance of the initial point cloud $Z_0$ as
follows:
\begin{equation}
h\coloneqq h_{Z_0,Q}\coloneqq 
\sup_{z\in Q}\inf_{z_i\in Z_0} \vert z - z_i\vert.
\end{equation}

The numerical approximation $f_{h,\sigma}(t,z)$ of $f(t,z)$ is then defined 
for all times $t\ge 0$ analogous to \eqref{eqn:kernel-system}, i.\,e.,
the time-dependent coefficients $\mathsf{c}(t)=\bigl(c_1(t),\dotsc,
c_N(t)\bigr)\in\reals^N$ are given as the solution of the linear system:
\begin{equation}
\label{eqn:coeffs_f_h_sigma}
\mathsf{K}(t)\mathsf{c}(t)=\mathsf{f}
\end{equation}
with time-dependent matrix entries $\mathsf{K}_{ij}(t) \coloneqq
k\bigl(z_i(t),z_j(t)\bigr)$ and
\begin{equation}
\label{eqn:sum_f_h_sigma}
f_{h,\sigma}(t,z) \coloneqq\sum_{i=1}^N c_i(t) k(z_i(t), z).
\end{equation}

We then obtain the following result.

\begin{theorem*}[Convergence linear case]
\label{thm:convergence_linear_case}
Let the above assumptions be fulfilled and let $T>0$. Let $m=n+1$. Then, for $h$ small emough
and for all $t\in[0,T]$, the interpolant $f_{h,\sigma}$  satisfies the error
bound
\begin{equation}
\label{eqn:error-bound-linear-case}
\Vert f(t,\cdot) - f_{h,\sigma}(t,\cdot) \Vert_{L^p(Q)} \le C(T)h^{m-l(p)} \Vert f_0 \Vert_{H^m(Q)},
\end{equation}
where $C(T)>0$ depends on the problem and on discretisation parameters
such as $\vmax$ and order $n$, but is \emph{independent of $h$}.
The constant $l\coloneqq l(p)$ is defined as
\begin{align*}
l(p) = \begin{cases}
1-\frac{2}{p} & \text{if }2\le p\le\infty, \\
0             & \text{if }1\le p\le 2.
\end{cases}
\end{align*}
\end{theorem*}
\begin{proof}
Let $t \ge 0$ be a fixed time and the exact solution of
\eqref{eqn:char1}--\eqref{eqn:char2} with the starting positions $Z(0)=Z_0$.
Furthermore assume that the particle positions $Z(t)$ are mapped back to their
periodically equivalent positions in $Q$ at time $t$, i.\,e., $Z(t) \subset Q$
without loss of generality. 

From \cref{section:vp-theory} we know that $f(t,\cdot)$ has the
same regularity as $f_0$. Thus, using standard estimates we obtain the
existence of a constant $C_1>0$ independent of $h$ such that:%
\autocite[Corollary~11.33]{wendland_scattered_2004}
\begin{equation}
\Vert f(t,\cdot)-f_{h,\sigma}(t,\cdot)\Vert_{L^{p}(Q)}
\le C_1 h_{Z(t),Q}^{m-l(p)}\Vert f(t,\cdot)\Vert_{H^{m}(Q)}.
\end{equation}
Using \eqref{eqn:function_values_f}, and H\"older's inequality one obtains
\begin{equation}
\Vert f(t,\cdot)\Vert_{H^m(Q)} = \Vert f_0\circ\Phi_t^0\Vert_{H^m(Q)}
\le \underbrace{\Vert\Phi_t^0\Vert_{H^m(Q)}}_{\eqqcolon C_2(t)}\Vert f_0\Vert_{H^m(Q)}.
\end{equation}
The constant $C_2(t)$ depends on time $t$ and the problem, but is independent of $h$
and thus:
\begin{equation}
\Vert f(t,\cdot) - f_{h,\sigma}(t,\cdot) \Vert_{L^p(Q)} 
\le C_1 C_2(t) h_{Z(t), Q}^{m - l(p)} \Vert f_0 \Vert_{H^m(Q)},
\end{equation}
where $h_{Z(t),Q}$ is the fill distance of $Z(t)$.\autocite[section 1]{cottet_particle_1984}

Finally, using that $\Vert\nabla\Phi_{0}^{t}\Vert_{L^\infty}$ is bounded for all
finite $t$, one obtains
\begin{equation}
h_{Z(t),\Omega}\le C_3(t)h,
\end{equation}
where the constant $C_3(t)>0$ again depends on the problem and time $t$.%
\autocite[Lemma~3-5]{cottet_particle_1984} Thus, letting:
\begin{equation*}
C(T)\coloneqq \sup_{t\in[0,T]} C_1C_2(t)C_3(t)^{m-l(p)}
\end{equation*}
the result follows. This concludes the proof.
\end{proof}

Our numerical experiments indicate that high orders of convergence are preserved in the 
non-linear setting. A complete proof is beyond the scope of this work.

\section{Numerical results}
We consider three standard benchmarks: 
\emph{weak linear Landau damping}, \emph{two stream instability} and 
\emph{bump on tail instability}. For the interpolation step we used both
the direct and piece-wise ansatz (PW), depending on the 
number of particles.
Piece-wise interpolants tend to develop overshoots near the boundaries of 
the respective boxes, but can be efficiently computed for large numbers of
particles. Direct kernel-interpolants use the entire set of particles at once
and therefore do not suffer from this problem. 
This, however, comes at the price of an $O(N^3)$ complexity and for this reason
the direct approach does not scale well to large numbers of particles.
We will carry out experiments using both approaches to assess whether the 
increased accuracy of the direct ansatz outweighs its cost. 

Tensorised Wendland kernels were used for all interpolations, i.\,e., we
used kernels of the shape:
\begin{equation*}
k(z,\tilde{z}) = b_{1,n}^{W}\left(\frac{|x-\tilde{x}|}{\sigma_x}\right)
                 b_{1,n}^{W}\left(\frac{|v-\tilde{v}|}{\sigma_v}\right).
\end{equation*}
In the following we will only specify the order $n\in\naturals$ of the kernel
in use.

Preliminary experiments have shown that for our test cases reasonable ranges 
for the scaling parameters are $0.5 \le \sigma_x, \sigma_v \le 6$.
For the finest discretisations the resulting kernel matrices became too 
ill-conditioned even for direct linear solvers. For this reason we apply 
Tikhonov regularisation and solve the modified systems $\left(\mathsf{K}
+\mu^2\mathsf{I}\right)\mathsf{c}=\mathsf{f}$ with regularisation parameter
$\mu = 10^{-6}$. 

For the phase-space sub-division we use $N_{\mathrm{min}}\in [100,200]$. The Poisson
solver is a standard Galerkin method using uses B-Splines of order~8 (degree~7) on a
uniform grid with $\Delta x = \frac{L}{256}$. The high resolution of the Poisson solver
was chosen such that we can neglect the influence of errors in the computation
of the electric field $E_{h,\sigma}$ from $\rho_{h,\sigma}$. Furthermore note that the
computation time of $E_{h,\sigma}$, even for this resolution, is neglible compared to the 
computation time of the interpolation step.

\subsection{Weak Landau damping}
\label{sec:weak_landau}

The initial condition is 
\begin{equation}
\label{eqn:lin-landau-f0}
f_0(x,v) := \frac{1}{\sqrt{2\pi}} e^{-\frac{v^2}{2}} \big{(} 1 + \alpha
\cos(kx)\big{)}, \quad (x,v) \in [0,L] \times \R
\end{equation}
with $k = 0.5$, $\alpha = 0.01$, $L = 4\pi$. The velocity space is cut at
$\vmax = 6$. The initial state \eqref{eqn:lin-landau-f0} 
is a small perturbation to the Maxwellian distribution 
\begin{align*}
f_M(v) = \frac{1}{\sqrt{2\pi}} e^{-\frac{v^2}{2}},
\end{align*}
which is a steady state solution of the Vlasov--Poisson equation \eqref{eqn:vlasov}. 
For the direct method we used
$h_x = \frac{L}{32}$ and $h_v=\frac{v_{\max}}{32}$,
while for the PW ansatz $h_x = \frac{L}{512}$ 
and $h_v = \frac{v_{\max}}{512}$ were used. The scaling parameters were chosen as 
$\sigma_x = 3$, $\sigma_v = 1$
for the direct method and $\sigma_x = 6$, $\sigma_v = 3$ for the PW method.
We used order $n=2$  for both methods. 
For time-integration we used the classical Runge--Kutta method with 
$\Delta t = \frac{1}{8}$ for the direct method and a symplectic Euler 
scheme with $\Delta t = \frac{1}{16}$ for the PW method. 
We chose a low order time-integration method for the 
PW ansatz as the reconstructed solution is only 
piecewise continuous. 
Note that while larger time-steps would be possible due to the lack of a
CFL-condition, we opted for smaller time-steps to better resolve 
the evolution of the electric field amplitude when plotting. 

\begin{figure}
\centering
\includegraphics{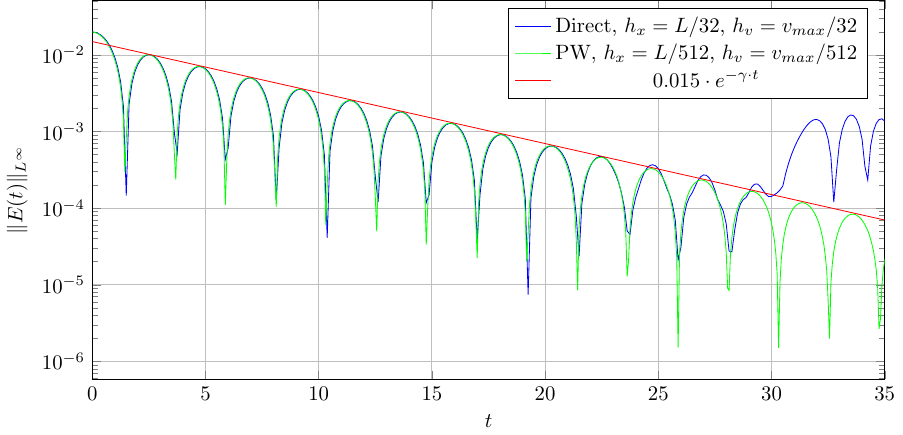}
\caption{\label{fig:lin_landau_amplitude}Electric field amplitude for the
\emph{weak Landau damping} benchmark. The direct method (blue graph)
uses $h_x = \frac{L}{32}$, $h_v = \frac{v_{\max}}{32}$, $\Delta t = \frac{1}{8}$, 
fourth order classical Runge--Kutta method and kernels of order $n=2$.
The PW method (green graph) uses $h_x = \frac{L}{512}$, 
$h_v = \frac{v_{\max}}{512}$, $\Delta t = \frac{1}{16}$, symplectic Euler
method and also kernels of order $n=2$. Up to $t= 25$ the numerical damping rate
of both methods are in good agreement  with the theoretical prediction. 
At $t \approx 29$ the recurrence effect sets in for the direct method, and 
between $t=24$ and $29$ the amplitude maxima are slightly overshooting their 
theoretical value. The PW method, on the other hand, reproduces the correct 
damping behaviour until $t = 35$.}
\end{figure}

\begin{figure}
\centering
\includegraphics{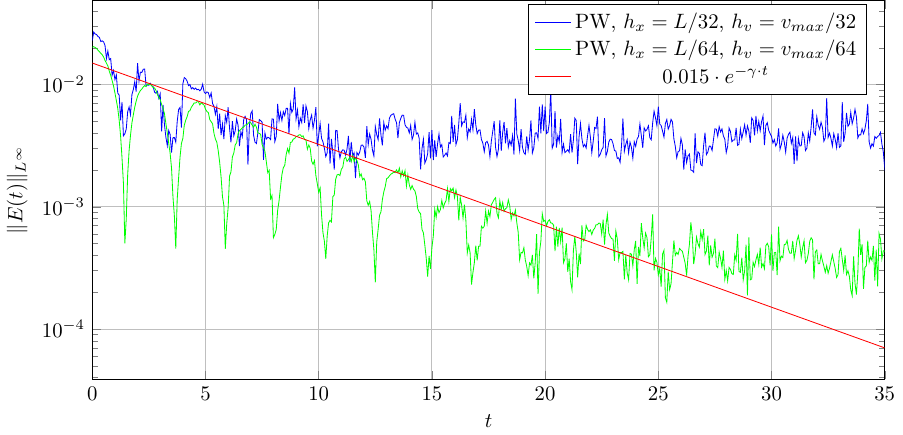}
\caption{\label{fig:e_amp_lin_landau_pw_with_low_res} Amplitude of $E_{h,\sigma}$ when running PW method 
with low resolution, resolution and parameters as for the direct method above.
}
\end{figure}

\Cref{fig:lin_landau_amplitude} shows the evolution of the electric field
amplitude over time. From theory and numerical experiments we know that the
electric field will get  damped periodically at rate $\gamma = 0.153359$ and
oscillation frequency  $\omega = 1.41566$\autocite{ARBER2002339}. 
Until $t\approx 25$ both methods are in good agreement with theory.
At this point we observe a recurrence phenomenon for the direct approach: 
the electric field amplitude increases again and a new damping process begins.
This is remarkable as recurrence has been known for grid-based 
solvers\autocite{mehrenberger:hal-01942708},
but is has been claimed that particle methods do not suffer from this 
effect\autocite{PhysRevE.84.036702}. 
Conventional blob methods, however, become \enquote{noisy} at this stage,
and we believe that it is this noise that masks the 
recurrence\autocite{myers_4th-order_2016,wang_particle--cell_2011,ameres_stochastic_2018}.

\begin{figure}
\centering
\begin{subfigure}{0.45\textwidth}
\includegraphics[width=\textwidth]{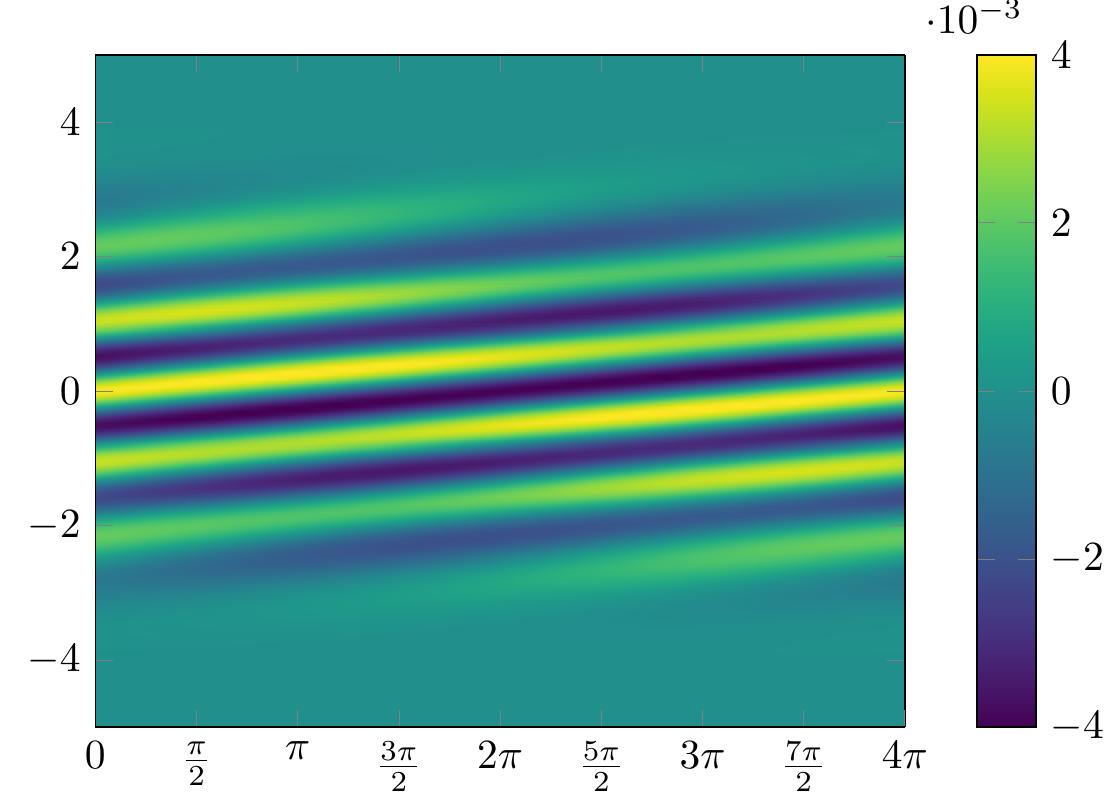}
\subcaption{Direct, $t=10$}
\label{subfig:lin_landau_direct_t_10}
\end{subfigure}
\begin{subfigure}{0.45\textwidth}
\includegraphics[width=\textwidth]{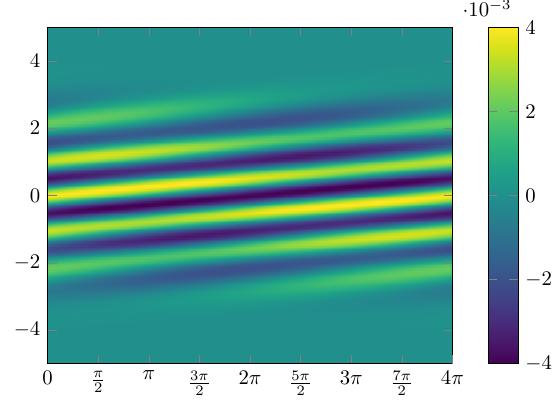}
\subcaption{PW, $t=10$}
\label{subfig:lin_landau_pw_t_10}
\end{subfigure}
\\[0.5cm]
\begin{subfigure}{0.45\textwidth}
\includegraphics[width=\textwidth]{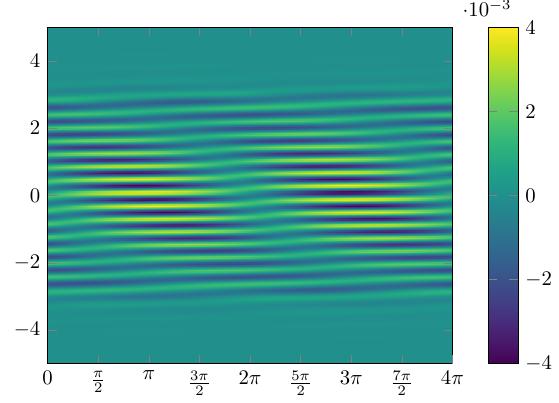}
\subcaption{Direct, $t=30$}
\label{subfig:lin_landau_direct_t_30}
\end{subfigure}
\begin{subfigure}{0.45\textwidth}
\includegraphics[width=\textwidth]{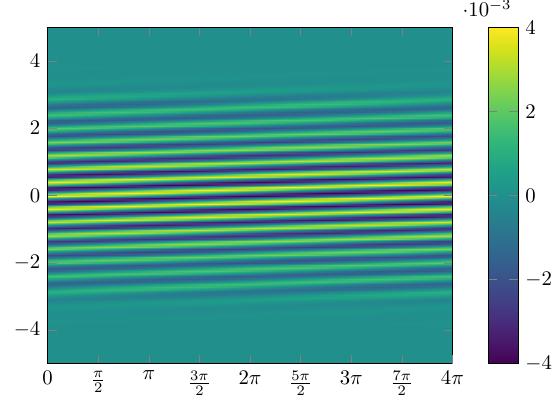}
\subcaption{PW, $t=30$}
\label{subfig:lin_landau_pw_t_30}
\end{subfigure}
\\[0.5cm]
\begin{subfigure}{0.45\textwidth}
\includegraphics[width=\textwidth]{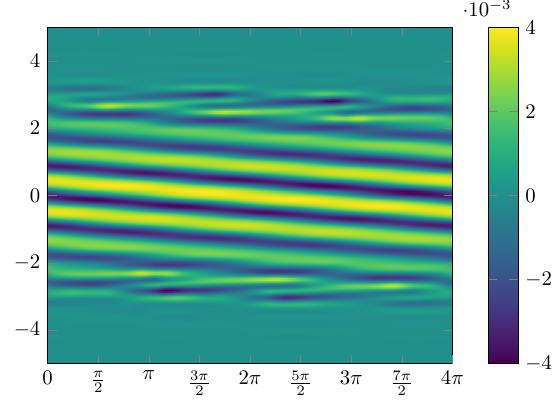}
\subcaption{Direct, $t=50$}
\label{subfig:lin_landau_direct_t_50}
\end{subfigure}
\begin{subfigure}{0.45\textwidth}
\includegraphics[width=\textwidth]{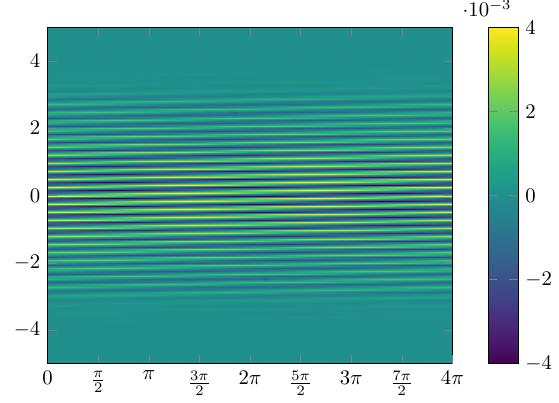}
\subcaption{PW, $t=50$}
\label{subfig:lin_landau_pw_t_50}
\end{subfigure}
\caption{Difference between the numerical solution $f_{h,\sigma}$ and the 
Maxwellian equilibrium $f_M$ 
for the \emph{weak Landau damping} benchmark. 
On the left-hand-side $f_{h,\sigma}$ is computed via the direct method and on 
the right-hand-side using the PW method. 
The resolutions are the same as in \cref{fig:lin_landau_amplitude}.
For $t \le 30$ the number of waves is steadily increasing. For the lower 
resolution, starting at $t=30$, numerical artefacts can be observed. At $t=50$ 
aliasing occurs: the spatial resolution is insufficient to capture the highest 
frequencies. Therefore high frequencies appear as low ones. 
The high-resolution solution can resolve the solution correctly up to $t=50$.}
\label{fig:linear_landau_f-maxwellian} 
\end{figure}

The cause of this phenomenon can be seen in
\cref{fig:linear_landau_f-maxwellian}, which shows the evolution of the
difference between the numerical solution $f_{h,\sigma}$ 
and the Maxwellian distribution $f_M$. One can observe waves of increasing 
frequency entering the domain from $v\rightarrow \pm\infty$. 
Starting at $t=30$, unphysical artefacts 
appear in the plot of the direct method. At this point we can observe
an aliasing effect: the high-frequency modes are not correctly captured
by the low resolution of the direct method.
The higher resolution of the PW method, however, can correctly
reproduce $f$ for extended period of times.

In the limit $t \rightarrow +\infty$ 
and for the initial datum \eqref{eqn:lin-landau-f0} the analytic solution 
$f(t,\cdot,\cdot)$  will weakly converge to a steady state but not in a strong 
sense\autocite{mouhot_landau_2011}. 
The solution $f$ develops waves of increasing frequency and number 
with time and therefore develops a small scale structure.
Thus the PW method will eventually suffer from recurrence as well. 

\begin{figure}
\centering
\begin{subfigure}{0.45\textwidth}
\includegraphics[width=\textwidth]{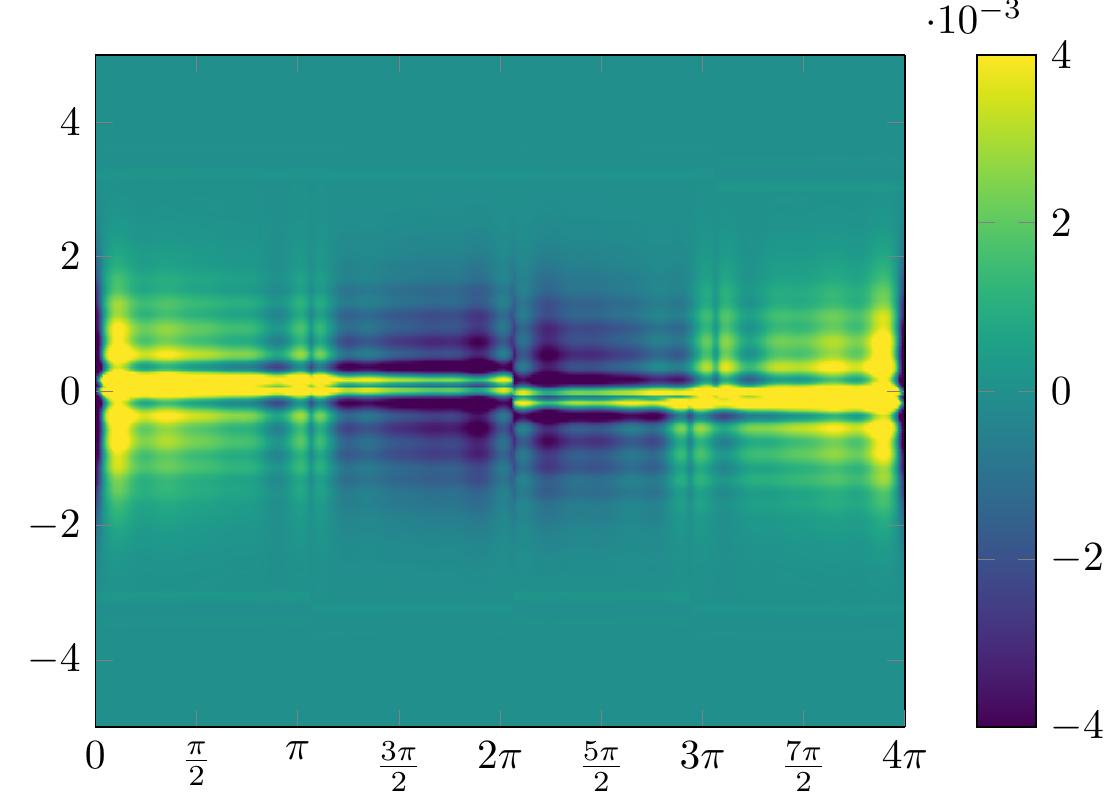}
\subcaption{$t=0$}
\label{subfig:lin_landau_pw_low_res_t_0}
\end{subfigure}
\begin{subfigure}{0.45\textwidth}
\includegraphics[width=\textwidth]{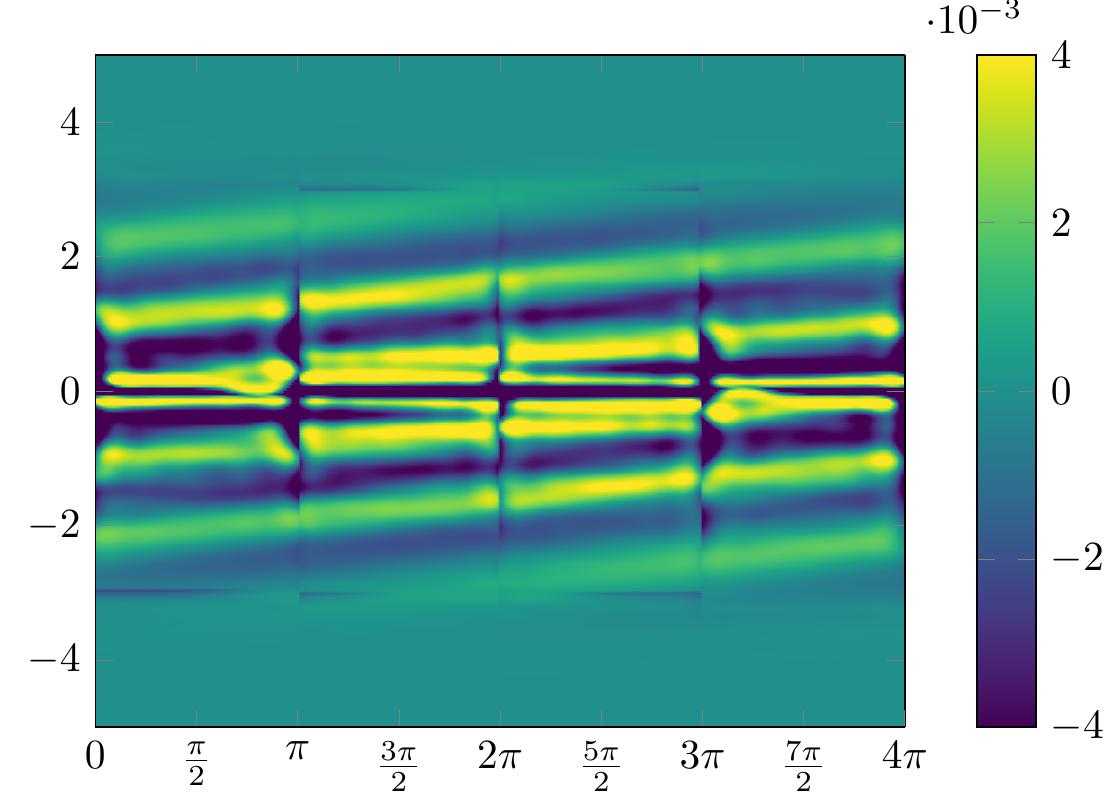}
\subcaption{$t=10$}
\label{subfig:lin_landau_pw_low_res_t_10}
\end{subfigure}
\\[0.5cm]
\begin{subfigure}{0.45\textwidth}
\includegraphics[width=\textwidth]{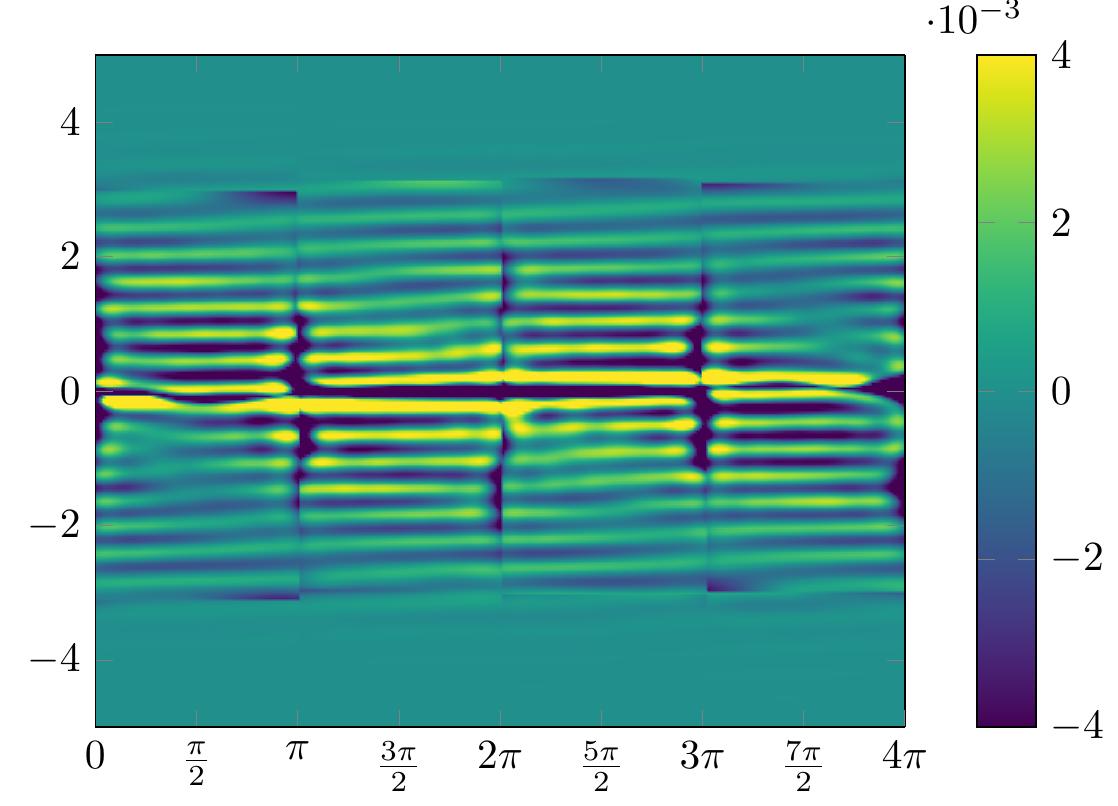}
\subcaption{$t=30$}
\label{subfig:lin_landau_pw_low_res_t_30}
\end{subfigure}
\begin{subfigure}{0.45\textwidth}
\includegraphics[width=\textwidth]{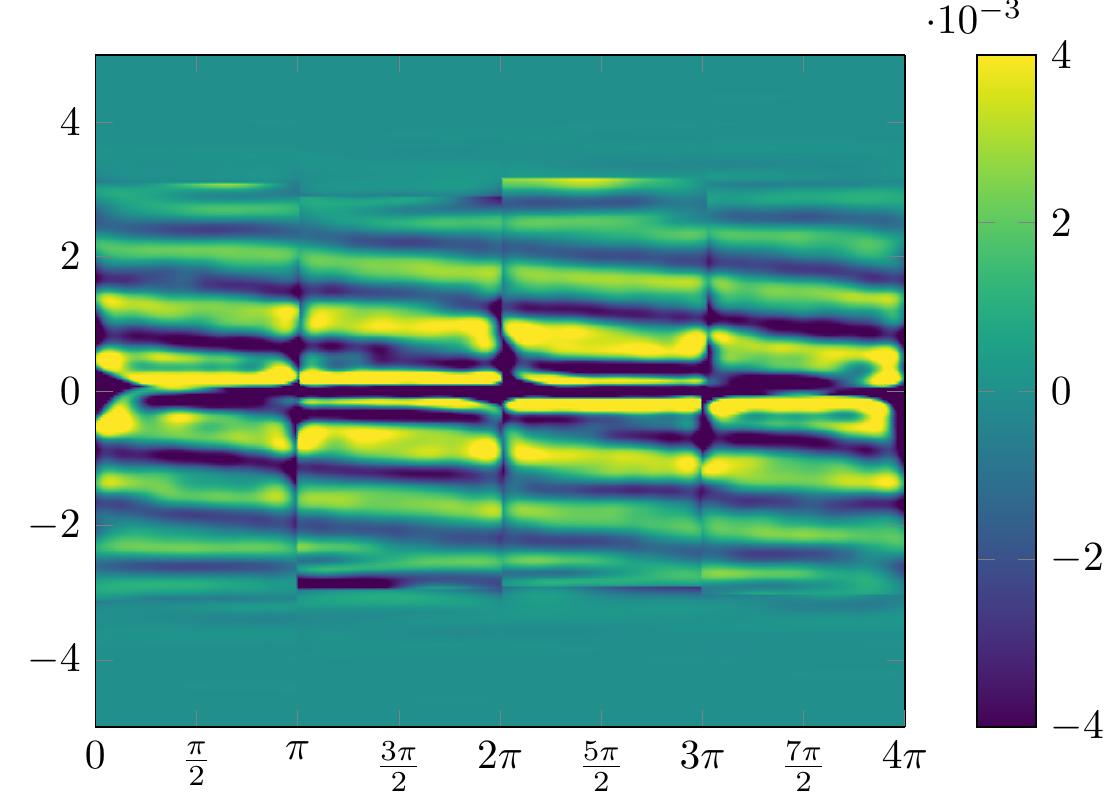}
\subcaption{$t=50$}
\label{subfig:lin_landau_pw_low_res_t_50}
\end{subfigure}
\caption{\label{fig:lin_landau_PW_with_low_res}Difference between $f_{h,\sigma}$ and $f_M$ for the
\textit{weak Landau damping} benchmark when using the PW method at a low resolution, equal to that
of the direct method. One observes jumps in the errors along the boundaries of the boxes $B$ used
in the sub-division of the phase space. However, the simulation still captures the correct global
dynamic of $f$. }
\end{figure}

In \Cref{fig:lin_landau_PW_with_low_res} we display the results of running the
PW method at low resolution; nameley the same number of particles and the same
scaling parameters as for the direct method. Along the boundaries of the boxes $B$
one observes errors caused by the discontinuity of the PW approximations. This
also leads to noisy results in the amplitude plot, \Cref{fig:e_amp_lin_landau_pw_with_low_res}.
Note however, that for a discontinuous function the $L^\infty$-norm is not suitable 
to  analyse errors and therefore the noisy results in below figure are to be expected
when using to low resolutions. Still even though there are strong errors along the
boundaries of the cover boxes, the overall dynamic in $f$ is correctly captured to
a similar extend as for the direct method. Especially one observes a similar
wave-structure for $t=10, 30$ and $50$ when comparing
\Cref{subfig:lin_landau_direct_t_10,subfig:lin_landau_direct_t_50} with 
\Cref{subfig:lin_landau_pw_low_res_t_10,subfig:lin_landau_pw_low_res_t_50}. 

This leads to the conclusion that using the PW method starts to be reasonable
when discontinuity errors are on the same level as the local interpolation errors,
which can only be expected for high enough resolutions. When looking at 
\Cref{fig:e_amp_lin_landau_pw_with_low_res} we see that doubling the resolution
from $h_x=L/32$ and $h_v=\vmax/32$ to $h_x=L/64$ and $h_v=\vmax/64$ already
significantly lowers the errors in the amplitude plots, suggesting that the latter
resolution is the lowest resolution to get reasonable results with PW methods for
this set of parameters and this benchmark. This also the reason why we chose to run
our tests for the PW method with $h_x=L/512$ and $h_v=\vmax/512$, i.\,e., a 
relatively high resolution. Note that the direct method can be run with lower
resolution and still produce good results, however, takes significantly longer
to run and uses significantly more resources. This is why we decided to run the 
test for the direct method in the relatively low resolution only. The performance 
and accuracy trade-off between the direct and PW method will be discussed in more detail
later in \Cref{subsec:conv_study} and \Cref{subsec:comp_efficiency}.

\begin{figure}
\centering
\begin{subfigure}{0.45\textwidth}
\includegraphics[width=\textwidth]{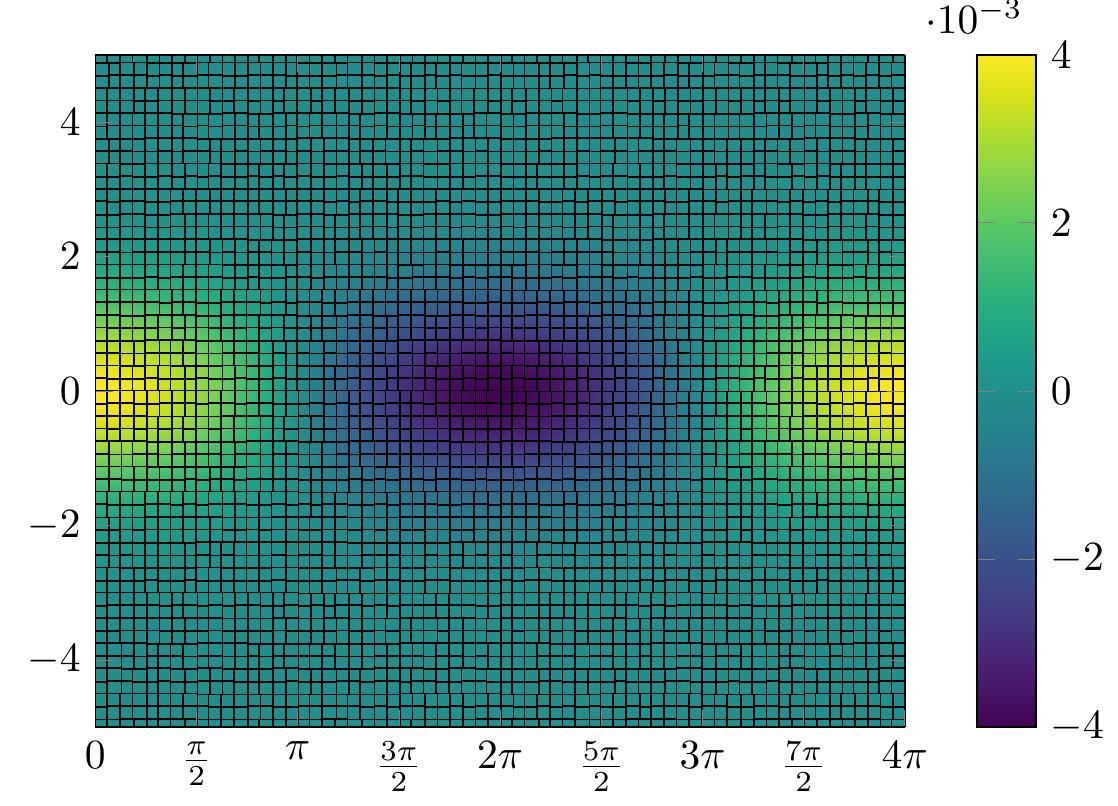}
\end{subfigure}
\begin{subfigure}{0.45\textwidth}
\includegraphics[width=\textwidth]{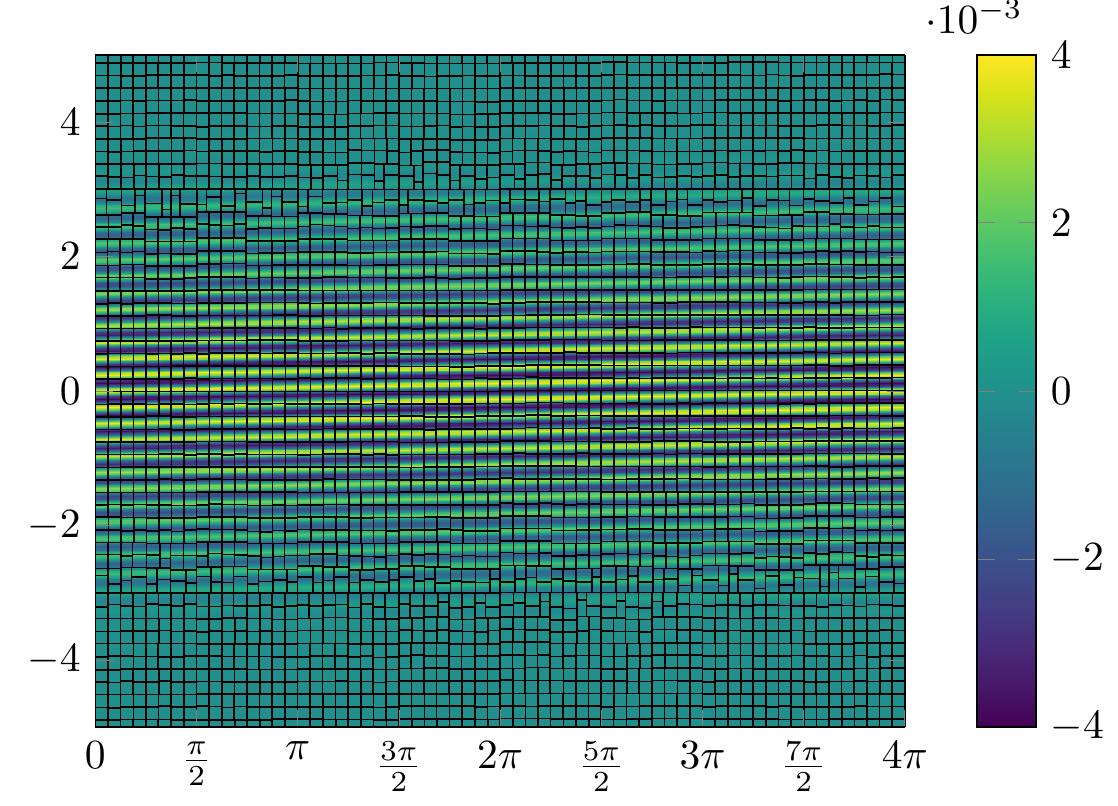}
\end{subfigure}
\caption{\label{fig:comparison_kd_tree_boxes_lin_landau_512_1024}
Comparison between domain decompositions for the PW method for $t=0$
and $t=50$ for the \textit{weak Landau damping} benchmark using 
$h_x = \frac{L}{512}$ and $h_x = \frac{\vmax}{512}$.}
\end{figure}

In figure \ref{fig:comparison_kd_tree_boxes_lin_landau_512_1024} we compare
the domain decomposition using the $kd$-tree for different times in the 
simulation using $h_x = \frac{L}{512}$ and $h_x = \frac{\vmax}{512}$.
We notice that the decomposition essentially stays close to an uniform grid 
for both $t=0$ and $t=50$. The strongest adaptation can be observed where waves 
enter the domain, around $\vert v \vert \approx 3$.

\subsection{Two stream instability}
For our second benchmark we consider the initial condition
\begin{equation}
\label{two-stream-f0}
f_0(x,v) = \frac{1}{\sqrt{2\pi}} v^2 e^{-\frac{v^2}{2}} \big{(}1 + \alpha
\cos(kx)\big{)}, \quad (x,v) \in [0,L] \times \R
\end{equation}
with $\alpha = 0.01$, $k = 0.5$, $L = 4\pi$. We use $v_{\max} = 8$ as cut-off in velocity-space.
The \emph{two stream instability} simulates two particle
streams of same density but opposing velocities colliding with each other. 
The initial state is a slight perturbation of the instable equilibrium
\begin{align*}
f_{eq}(v) = \frac{1}{\sqrt{2\pi}} v^2 e^{-\frac{v^2}{2}}.
\end{align*}

We used $h_x= \frac{L}{64}$ and $h_v = \frac{v_{\max}}{128}$ for the direct method
and $h_x= \frac{L}{512}$ and $h_v = \frac{v_{\max}}{1024}$ for the PW method. 
The scaling parameters were fixed as $\sigma_x = 2$, $\sigma_v = 1$
for the direct method and $\sigma_x = 4$, $\sigma_v = 2$ for the PW method.
For both methods we used fourth order kernels.
The time-steps are set to $\Delta t = \frac{1}{4}$ and $\Delta t = \frac{1}{32}$ 
with classical Runge-Kutta and symplectic Euler as time-integrators respectively.

\begin{figure}
\centering
\includegraphics{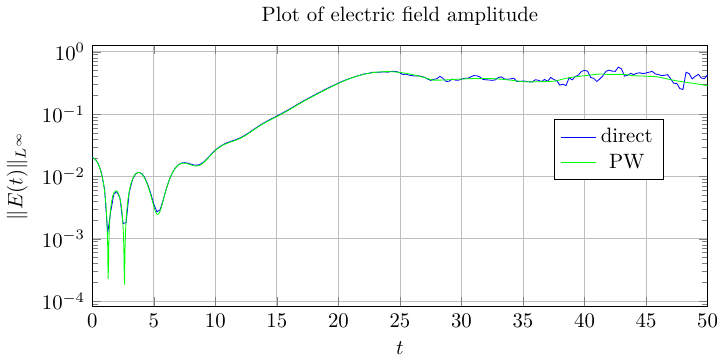}
\caption{\label{fig:two_stream_E}Electric field amplitude for the \emph{two
stream instability}. The direct method (blue) used $h_x = \frac{L}{64}$,
$h_v=\frac{10}{128}$, $\Delta t = \frac{1}{4}$ and the classical Runge-Kutta
method. The PW method (green) used $h_x = \frac{L}{512}$,  $\frac{10}{1024}$,
$\Delta t = \frac{1}{32}$ and the symplectic Euler method. Both the direct and
PW method employed fourth order kernels.  After an initial damping until $t
\approx 5$ the mixing process starts dominating the  dynamics. The electric
field amplitude arrives at its maximum at approximately  $t \approx 23$ and
starts periodically oscillating afterwards.  Until $t\approx 23$ the increase
of amplitude is captured correctly  by both methods.  After $t \approx 25$
slight numerical artefacts appear in the solution computed with lower
resolution.}
\end{figure}

The evolution of the electric field's amplitude is depicted in
\cref{fig:two_stream_E}. Because the direct approach is limited to fairly 
small numbers of particles, i.\,e., coarse resolutions, it performs 
worse than the PW method.
However, the direct method is still able to capture the dynamics
of both $f$ and $E$ qualitatively correctly.
The forming of the filaments introduces steep 
gradients. This results in overshoots in the numerical solution, as soon as
these gradients can no longer be resolved by the fixed resolution. 

\begin{figure}
\centering
\begin{subfigure}{0.45\textwidth}
\includegraphics[width=\textwidth]{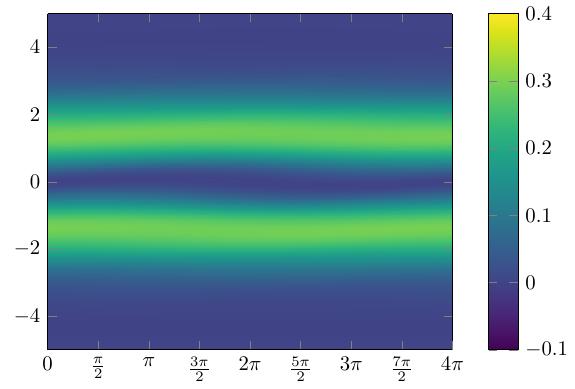}
\subcaption*{Direct, $t=10$}
\end{subfigure}
\begin{subfigure}{0.45\textwidth}
\includegraphics[width=\textwidth]{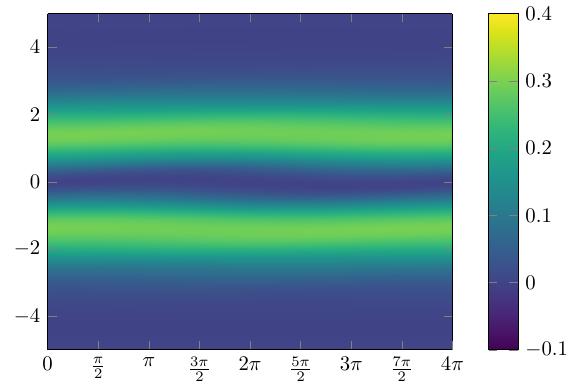}
\subcaption*{PW, $t=10$}
\end{subfigure}
\\[0.5cm]
\begin{subfigure}{0.45\textwidth}
\includegraphics[width=\textwidth]{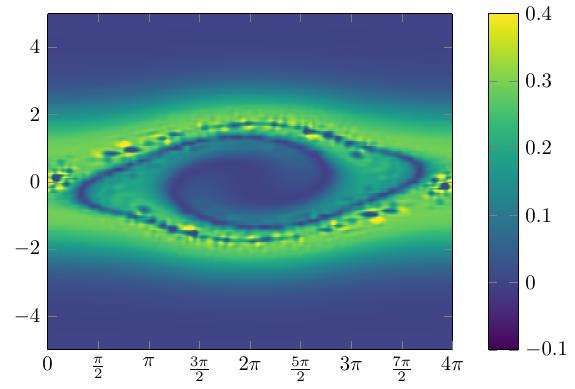}
\subcaption*{Direct, $t=30$}
\end{subfigure}
\begin{subfigure}{0.45\textwidth}
\includegraphics[width=\textwidth]{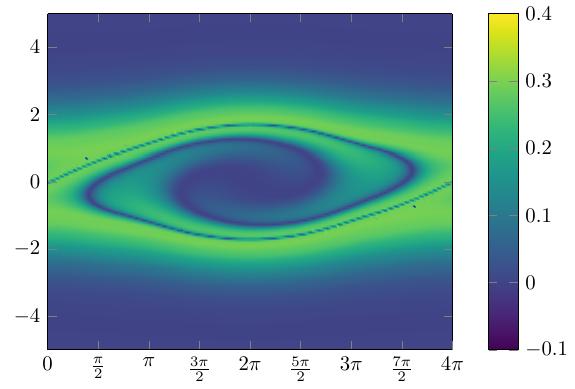}
\subcaption*{PW, $t=30$}
\end{subfigure}
\\[0.5cm]
\begin{subfigure}{0.45\textwidth}
\includegraphics[width=\textwidth]{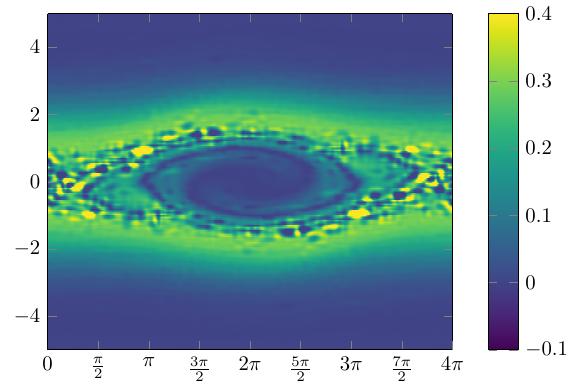}
\subcaption*{Direct, $t=50$}
\end{subfigure}
\begin{subfigure}{0.45\textwidth}
\includegraphics[width=\textwidth]{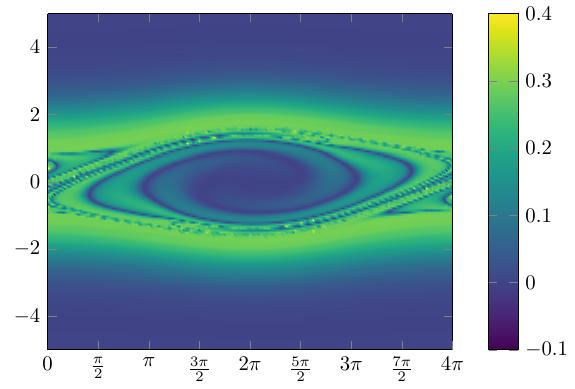}
\subcaption*{PW, $t=50$}
\end{subfigure}
\\[0.5cm]
\caption{\label{fig:two_stream_f}The distribution function
$f_{h,\sigma}(t,x,v)$ for the \emph{two stream  instability benchmark}. The
numerical solution computed using the direct method is  on the left-hand side,
the PW method on the right-hand side. The parameters were chosen as in
\cref{fig:two_stream_E}.  At $t \approx 10 $ the phase-space \enquote{vortex}
starts to form.  After $t=30$ one can observe the \enquote{vortex} rotation. 
With each rotation further filaments enter the \enquote{vortex}. It can also be
observed that  the filaments enter at $x=0$ or $x=4\pi$ and drift towards the
center of the \enquote{vortex}. While both methods reproduce the dynamics of
$f$ qualitatively correct, the errors made by the  direct method are
significantly higher. Both methods suffer from overshoots near filaments.} 
\end{figure}

This can cleary be observed in \cref{fig:two_stream_f} after $t \approx 30$  
in the plots for the direct method. While the direct and PW methods, similar to
grid-based methods, suffer from overshoots and thus do not provide good
accuracy for extended times in turbulent simulations,  they reproduce the fine
details of $f$ better than Eulerian or PIC methods. 

\begin{figure}
\centering
\begin{subfigure}{0.45\textwidth}
\includegraphics[width=\textwidth]{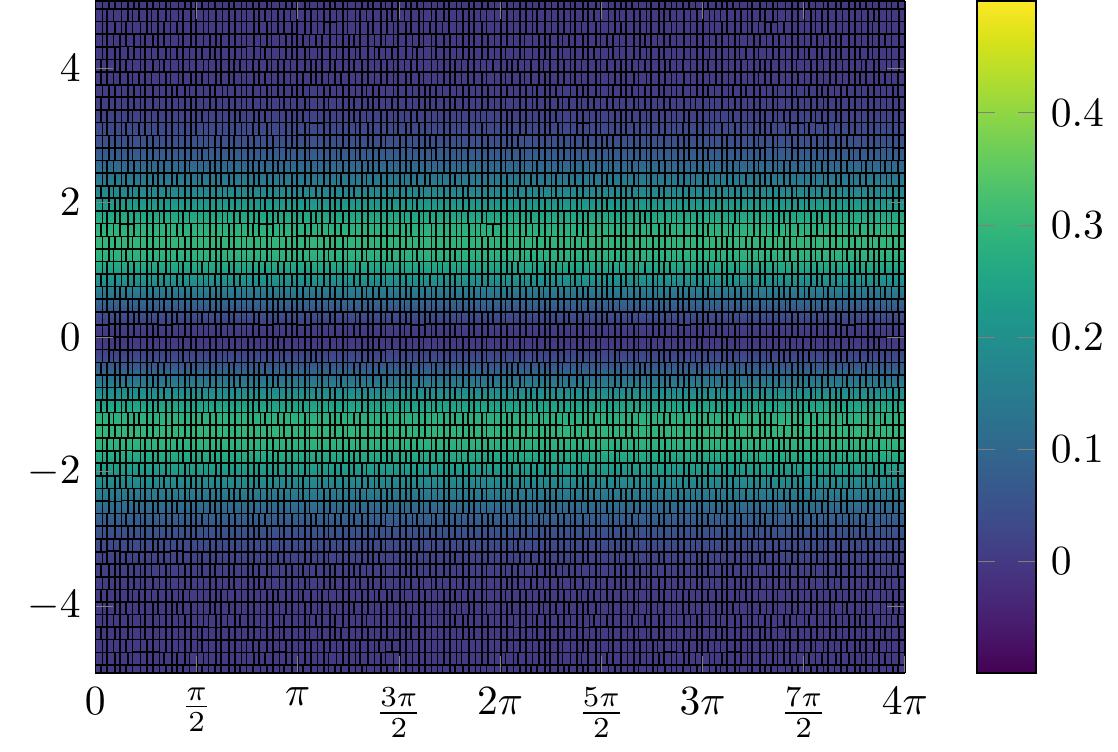}
\end{subfigure}
\begin{subfigure}{0.45\textwidth}
\includegraphics[width=\textwidth]{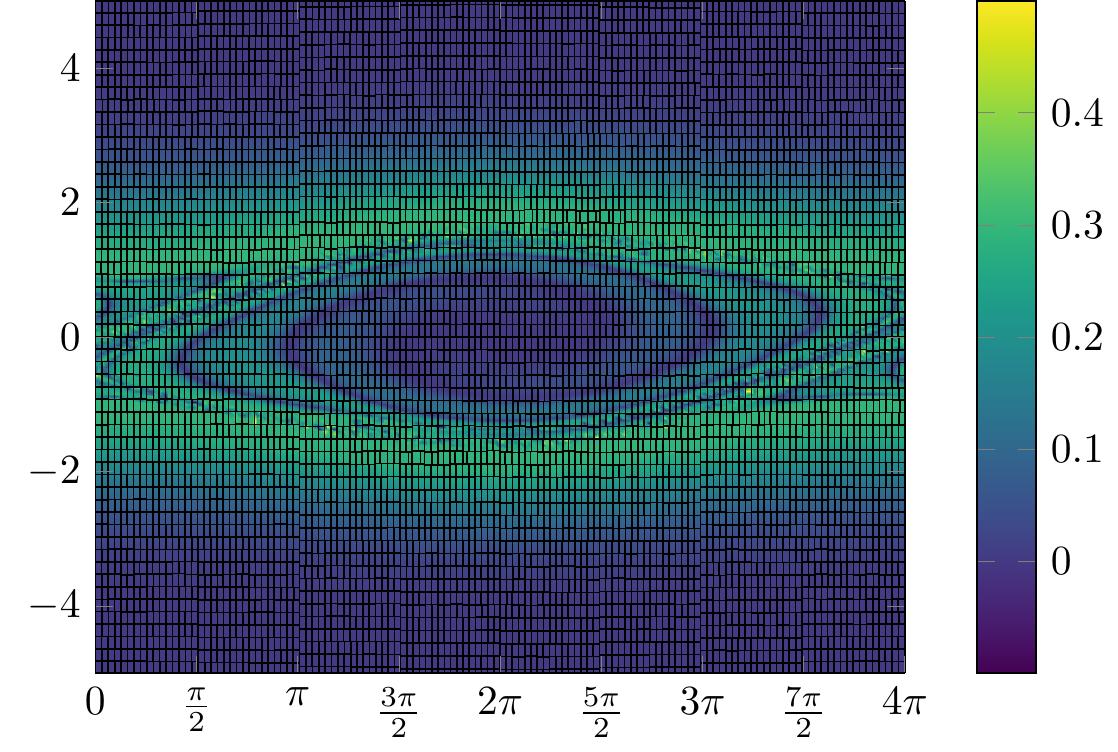}
\end{subfigure}
\caption{\label{fig:comparison_kd_tree_boxes_two_stream_1024_1024}
Comparison between domain decompositions for the PW method for $t=0$
and $t=50$ for the \textit{two stream} benchmark using 
$h_x = \frac{L}{512}$ and $h_x = \frac{\vmax}{1024}$.}
\end{figure}

In \Cref{fig:comparison_kd_tree_boxes_two_stream_1024_1024} we compare
the domain decompositions of the PW scheme for different times $t$ using
$h_x = \frac{L}{512}$ and $h_x = \frac{\vmax}{1024}$. At $t=0$ it is close to
a uniform grid. This is expected due to the uniforml particle distribution at
$t=0$. Later, at $t=50$ the particles are in more disarray and thus we observe
some adaptation of the domain decomposition. However, the decompositions seem
to be still close to an uniform grid, suggesting that the particle distribution
is also still quasi-uniform. The strongest adaptations can be seen
close to the filaments.

\subsection{Bump on tail instability}
The initial condition for our final benchmark, the \emph{bump on tail instability}, is
\begin{equation}
f_0(x,v) = \frac{1}{\sqrt{2\pi}} \left( n_p \cdot e^{-\frac{v^2}{2}} 
+ n_b \cdot e^{-\frac{1}{2}\frac{(v-v_b)^2}{v_t^2}}\right) \left( 1 + \alpha \cos(kx)\right).
\end{equation}
with $n_p = 0.9$, $n_b = 0.2$, $v_b = 4.5$, $v_t = 0.5$, $\alpha = 0.04$, $k = 0.3$,
$L = \frac{2\pi}{0.3}$\autocite{ARBER2002339}. The cut-off in velocity space is set to 
$v_{\max} = 10$. For this test case we only consider the PW method and are interested 
in both long- and short-term accuracy.
We chose $h_x = \frac{L}{1024}$ and $h_v = \frac{v_{\max}}{512}$ 
as resolution and second order kernels. The scaling parameters were set to
$\sigma_x = 6$, $\sigma_v = 3$. The time-step is set to $\Delta t = \frac{1}{16}$ 
and the symplectic Euler time integration method was used.

The \emph{bump on tail instability} simulates the clash of a low density particle stream
with Maxwellian velocity distribution around $v_b = 4.5$ into resting particles, i.\,e., 
with Maxwellian velocity distribution around $0$. The equilibrium state
\begin{align*}
f_{eq}(v) = \frac{1}{\sqrt{2\pi}} \left( n_p \cdot e^{-\frac{v^2}{2}} 
+ n_b \cdot e^{-\frac{1}{2}\frac{(v-v_b)^2}{v_t^2}}\right)
\end{align*}
gets slightly perturbed, whereby the mixing process is initiated. 
The resulting dynamics can be described as an overlapping of the effects of
\emph{weak Landau damping} and the \emph{two stream instability}, 
i.\,e., an overlapping of mixing and damping processes. 
Which effect dominates, depends on the difference in density and the strength 
of the initial perturbation. 

\begin{figure}
\centering
\includegraphics{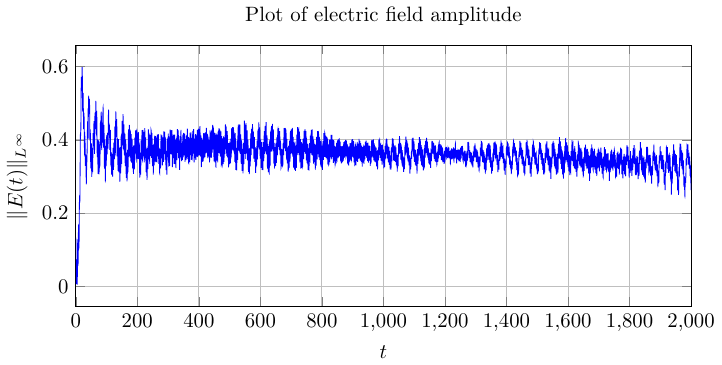}
\caption{\label{fig:bump_on_tail_piecewise_pou_long_times}Amplitude of the
electric field until $T= 2000$ for the \emph{bump  on tail instability}
benchmark. The simulation was run at  resolution $h_x = \frac{L}{1024}$, $h_v =
\frac{10}{512}$ and time step $\Delta t = \frac{1}{16}$ using the symplectic
Euler method.  After an initial increase of amplitude, several oscillation
modes can be observed.  In the long time limit $t \rightarrow 2000$ the
amplitude gets gradually damped.}
\vspace{1cm}
\includegraphics{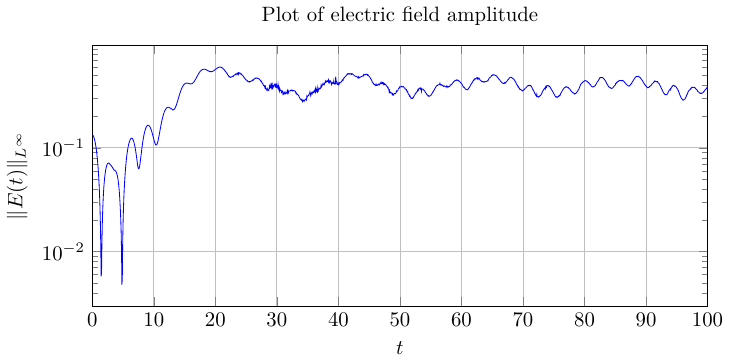}
\caption{\label{fig:bump_on_tail_E_pou_short}Amplitude of the electric field of
the \emph{bump on tail instability} until $T = 50$ with simulation parameters
as in figure \cref{fig:bump_on_tail_piecewise_pou_long_times}. After an initial
damping until  $t\approx 5$ the amplitude increases to reach its global maximum
at $t\approx 20$. Then there are two dominating oscillation modes, a slow
oscillation with period $\approx 22$ and another faster oscillation with period
$\approx 2.5$. The first oscillation is caused by the mixing of the two
particle streams, see \cref{fig:two_stream_E}, the second oscillation
is caused by a Landau damping effect, see \cref{fig:lin_landau_amplitude}. 
Between $t = 30$ and $t = 40$ one can observe numerical noise. At this point
the mixing of the two streams is causing steep gradients and thus overshoots
in the numerical solution $f_{h,\sigma}$ which results in numerical errors
when computing $E_{h,\sigma}$. After $t=50$, the amplitude is reproduced
correctly again.}
\end{figure}

For the chosen set of parameters, the damping effect on electric field is small
and  only becomes apparent at large time intervals. Therefore the benchmark
involves  being able to simulate the correct behaviour for times $t \gg 100$. In
\cref{fig:bump_on_tail_piecewise_pou_long_times} one sees the long time
damping,  which is in good agreement with the results presented by Arber and
Vann\autocite{ARBER2002339}.

\begin{figure}
\centering
\begin{subfigure}{0.45\textwidth}
\includegraphics[width=\textwidth]{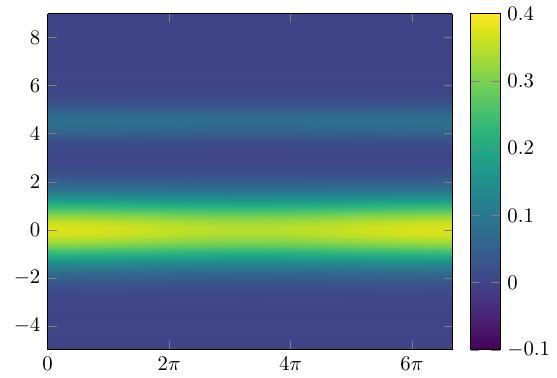}
\subcaption{$t=0$}
\label{subfig:bump_t_0}
\end{subfigure}
\begin{subfigure}{0.45\textwidth}
\includegraphics[width=\textwidth]{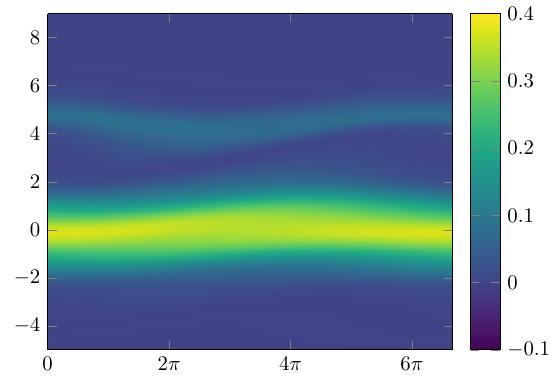}
\subcaption{$t=10$}
\label{subfig:bump_t_10}
\end{subfigure}
\\[0.5cm]
\begin{subfigure}{0.45\textwidth}
\includegraphics[width=\textwidth]{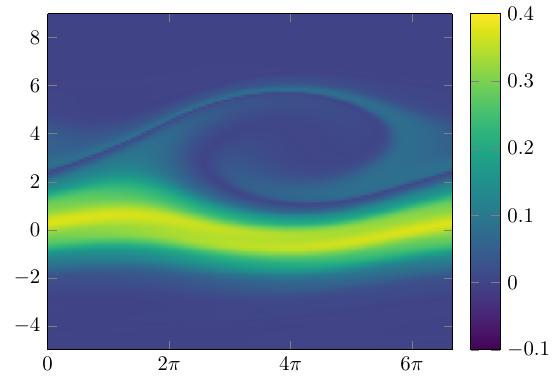}
\subcaption{$t=20$}
\label{subfig:bump_t_20}
\end{subfigure}
\begin{subfigure}{0.45\textwidth}
\includegraphics[width=\textwidth]{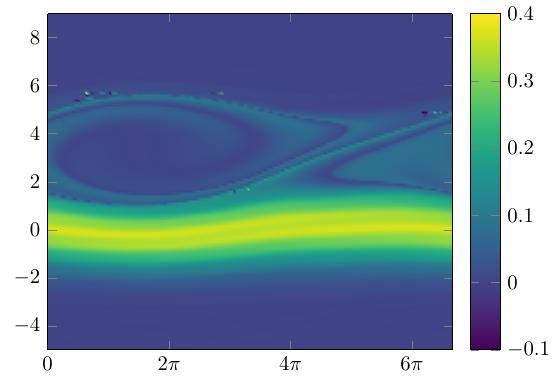}
\subcaption{$t=30$}
\label{subfig:bump_t_30}
\end{subfigure}
\\[0.5cm]
\begin{subfigure}{0.45\textwidth}
\includegraphics[width=\textwidth]{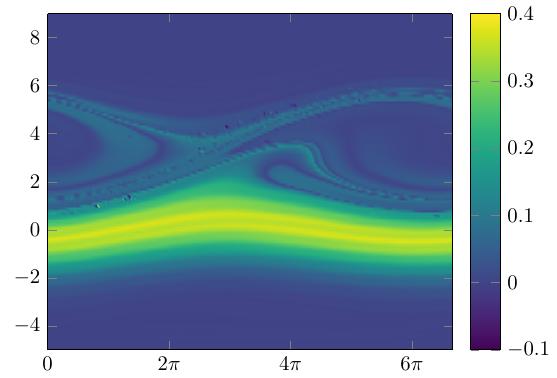}
\subcaption{$t=40$}
\label{subfig:bump_t_40}
\end{subfigure}
\begin{subfigure}{0.45\textwidth}
\includegraphics[width=\textwidth]{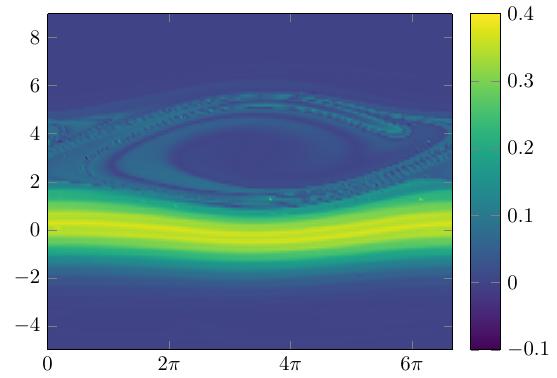}
\subcaption{$t=50$}
\label{subfig:bump_t_50}
\end{subfigure}
\caption{\label{fig:bump_on_tail_f_piecewise_pou}Evolution of the distribution
function $f_{h,\sigma}(t,x,v)$ until $T=50$ for the \emph{bump on tail
instability}. Computed using the PW method with the same parameters as in
\cref{fig:bump_on_tail_piecewise_pou_long_times}. Beginning at $t\approx 10$
the phase-space \enquote{vortex} starts to form and gets fully developed at $t
\approx 20$. Afterwards it starts periodically moving along the x-axis with a
period of $\approx 22$ which coincides with  the slow oscillation in
\cref{fig:bump_on_tail_E_pou_short}.  After $t=20$ the formation of so-called
Langmuir waves on the bigger particle  cluster can be observed. This coincides
with the faster oscillations appearing in \cref{fig:bump_on_tail_E_pou_short}
with period $\approx 2.5$.}
\end{figure}

In \cref{fig:bump_on_tail_E_pou_short} we see an initial damping between $t=0$
and $t \approx 5$. On the one hand, after $t \approx 5$ a \enquote{vortex}
begins to form that will eventually  dominate the dynamics for early times. The
\enquote{vortex} moves periodically in phase space along the  $x$-axis. This
can be seen in \cref{fig:bump_on_tail_f_piecewise_pou}.  After $t \approx 20$
an increase in the number of filaments can be observed.  This effect is similar
to that of \cref{fig:two_stream_f} from the \emph{two stream instability}
simulation. 

On the other hand, starting at $t \approx 20$ one can observe waves forming 
on the particle cluster centred at $v=0$. This is comparable to \emph{weak
linear  Landau damping}, see \cref{fig:linear_landau_f-maxwellian}. Note, that
compared  to the previous two benchmarks, the perturbation strength $\alpha =
0.04$ is  significantly higher and thus the amplitudes of appearing waves are
higher as well. 

Similar to the \emph{two stream instability} benchmark, we observe overshoots
resulting in slight numerical noise in this simulation, see
\cref{subfig:bump_t_40,subfig:bump_t_50}. But in contrast to the \emph{two
stream instability} benchmark, the gradients do not get as steep and therefore
the simulation stays stable even with fewer particles and for extended periods
of time.

\subsection{\label{subsec:conv_study}Convergence study}
In this section we investigate the convergence behaviour of the presented
methods in the  setting of the \emph{weak Landau damping} benchmark, see
\Cref{sec:weak_landau}. In particular we compare both the performance and 
accuracy of the direct and PW methods. We compute the errors for
$f_{h,\sigma}$ with respect to a reference solution with $h_x=\tfrac{L}{4096}$,
$h_v = \tfrac{v_{\mathrm{max}}}{4096}$ using the PW method and second order
kernels. Note that we use the PW method for the reference solution as computation
of a global interpolant in this resolution is not  feasible from both memory usage
and run-time perspective on the machine we used.

For all benchmarks we used the classical Runge--Kutta method with
$\Delta t = \frac{1}{8}$ for the direct approach and the symplectic Euler 
time-integration scheme with $\Delta t = \frac{1}{16}$  for the PW approach.
The minimal box-size for the PW method was set to $N_{\mathrm{min}} = 100$. 

\begin{figure}
\centering
\begin{subfigure}{0.48\textwidth}
\includegraphics[width=0.98\textwidth]{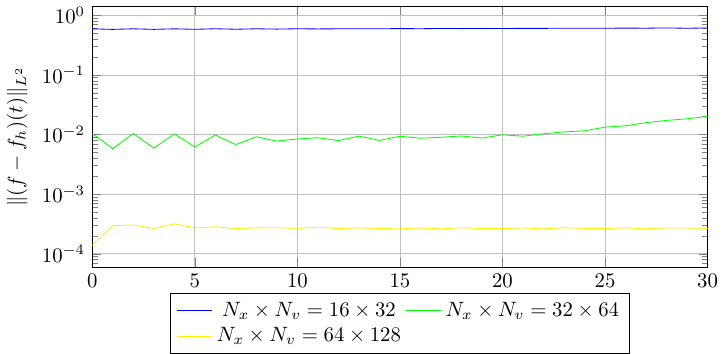}
\vspace{3mm}
\subcaption{\label{subfig:f_l2_error_2_order_kernel_direct} Direct}
\end{subfigure}
\begin{subfigure}{0.48\textwidth}
\includegraphics[width=0.98\textwidth]{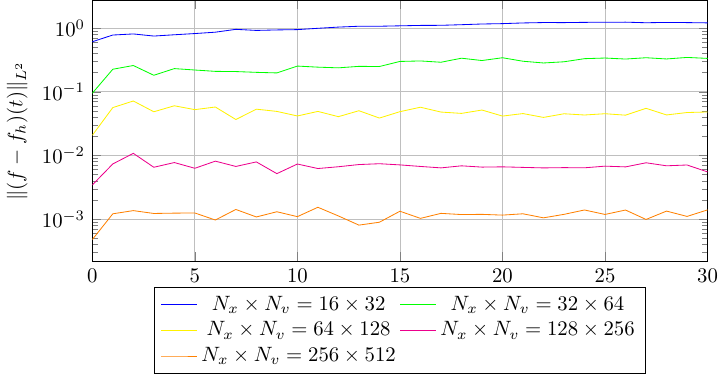}
\subcaption{\label{subfig:f_l2_error_2_order_kernel_pw} PW}
\end{subfigure}
\caption{\label{fig:f_l2_error_2_order_kernel}
A comparison of $L^2$-errors for $f_{h,\sigma}$ in time using second order 
kernels and $\sigma_x = 1$ and $\sigma_v = 0.5$. 
Left the direct and on the right the PW approach was used.
One observes a higher convergence order for the direct approach exceeding 
the theoretically predicted convergence order of $O(h^3)$, while the observed convergence order
for the PW method is slightly below the predicted convergence order.} 
\end{figure}

\begin{figure}
\centering
\begin{subfigure}{0.48\textwidth}
\includegraphics[width=0.98\textwidth]{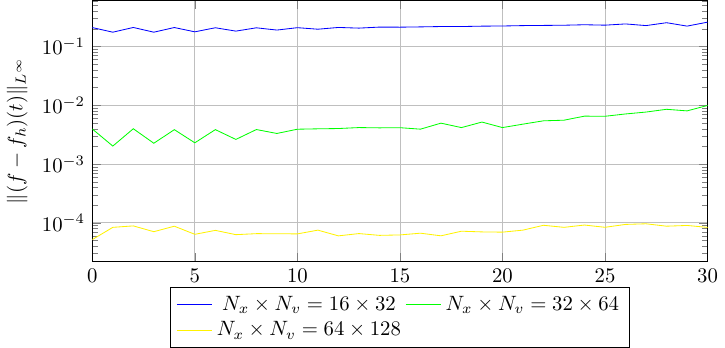}
\vspace{3mm}
\subcaption{\label{subfig:f_linfty_error_2_order_kernel_direct} Direct}
\end{subfigure}
\begin{subfigure}{0.48\textwidth}
\includegraphics[width=0.98\textwidth]{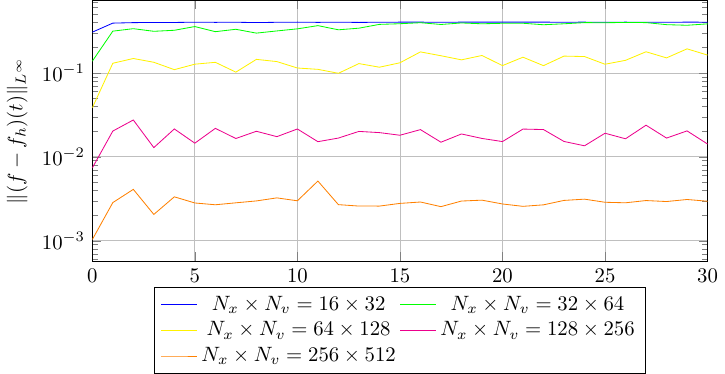}
\subcaption{\label{subfig:f_linfty_error_2_order_kernel_pw} PW}
\end{subfigure}
\caption{\label{fig:f_linfty_error_2_order_kernel}
A comparison of $L^\infty$-errors for $f_{h,\sigma}$ in time using second order kernels and $\sigma_x = 1$ and $\sigma_v = 0.5$.
Left the direct and on the right the PW approach was used.
One observes a higher convergence order for the direct approach, however, both approaches
exceed the theoretically predicted convergence order of $O(h^{2.5})$.} 
\end{figure}

\begin{figure}
\centering
\begin{subfigure}{0.48\textwidth}
\includegraphics[width=0.98\textwidth]{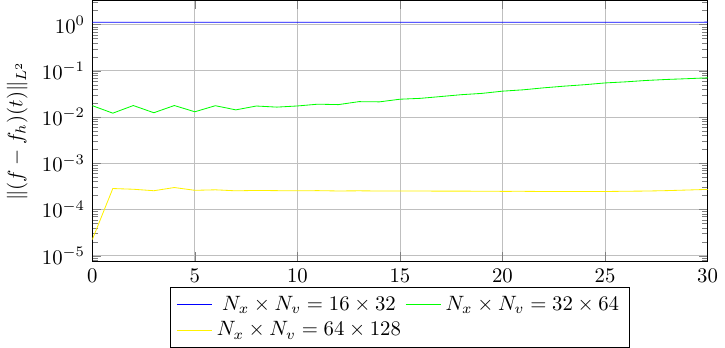}
\vspace{3mm}
\subcaption{\label{subfig:f_l2_error_4_order_kernel_direct} Direct}
\end{subfigure}
\begin{subfigure}{0.48\textwidth}
\includegraphics[width=0.98\textwidth]{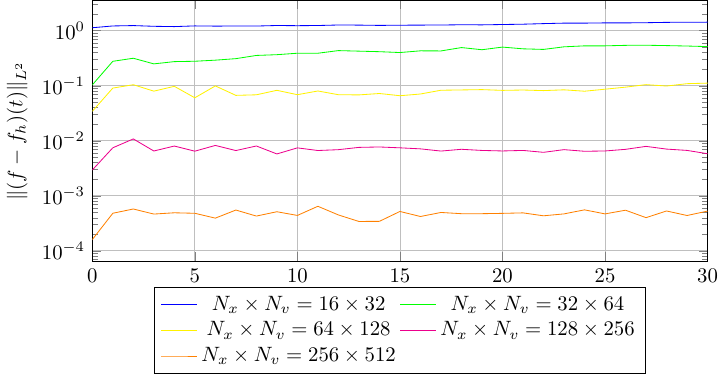}
\subcaption{\label{subfig:f_l2_error_4_order_kernel_pw} PW}
\end{subfigure}
\caption{\label{fig:f_l2_error_4_order_kernel}
A comparison of $L^2$-errors for $f_{h,\sigma}$ in time using fourth order kernels and $\sigma_x = 1$ and $\sigma_v = 0.5$.
Left the direct and on the right the PW approach was used.
One observes a higher convergence order for the direct approach again
exceeding the theoretically predicted convergence order of $O(h^5)$. 
For the PW method one observes a slighlty lower convergence order than 
expected.} 
\end{figure}

\begin{figure}
\centering
\begin{subfigure}{0.48\textwidth}
\includegraphics[width=0.98\textwidth]{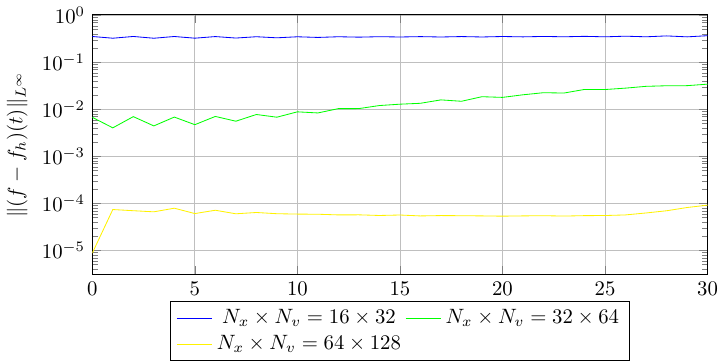}
\vspace{3mm}
\subcaption{\label{subfig:f_linfty_error_4_order_kernel_direct} Direct}
\end{subfigure}
\begin{subfigure}{0.48\textwidth}
\includegraphics[width=0.98\textwidth]{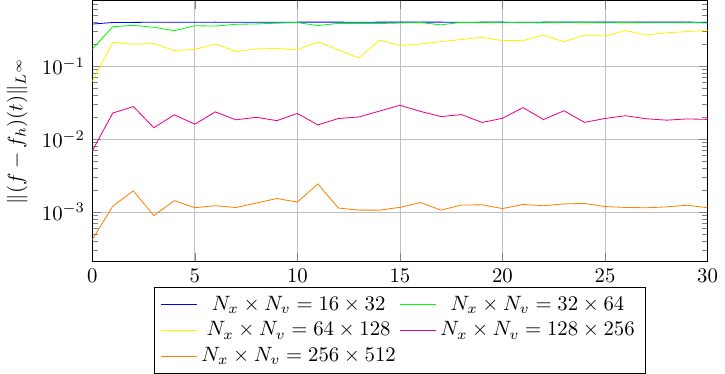}
\subcaption{\label{subfig:f_linfty_error_4_order_kernel_pw} PW}
\end{subfigure}
\caption{\label{fig:f_linfty_error_4_order_kernel}
A comparison of $L^\infty$-errors for $f_{h,\sigma}$ in time using fourth order kernels and $\sigma_x = 1$ and $\sigma_v = 0.5$.
Left the direct and on the right the PW approach was used.
One observes a higher convergence order for the direct approach 
exceeding the theoretically predicted convergence order of $O(h^{4.5})$. 
For the PW method one observes a slighlty lower convergence order than 
expected.} 
\end{figure}

The $L^2$- and $L^\infty$-errors for both schemes when using kernels of
order $n=2$ are displayed in
\Cref{fig:f_l2_error_2_order_kernel,fig:f_linfty_error_2_order_kernel}.
The corresponding results for kernels of order $n=4$ are given in
\Cref{fig:f_l2_error_4_order_kernel,fig:f_linfty_error_4_order_kernel}.

For all of these approaches, the error does not significantly grow
with time $t$; no \enquote{noise} becomes visibible. For the direct
method the $L^2$ and $L^\infty$-errors  are of similar magnitude.
Surprisingly, this is also the case for the $PW$ method, albeit to a
somewhat lesser extent. This lesser extent can be attributed to the fact
that the discontinuities of the PW approach more strongly affect the 
$L^\infty$-norm than the $L^2$-norm.

For the direct method we observe a rapid decrase of the error with
increasing resoltion. We empirically observe convergence rates of
$\mathcal{O}(h^{5.5})$ and $\mathcal{O}(h^7)$ in the $L^2$-norm
for respectively $n=2$ and $n=4$, clearly exceeding the expected rates.

For the PW approach, convergence only starts later, at higher resolutions
and does not reach the same rates as the direct method does. Here empirical
convergence rates approach $\mathcal{O}(h^{3})$ and $\mathcal{O}(h^{4.5})$ in 
the $L^2$-norm respectively for orders $n=2$ and $n=4$; much closer to the expected 
rates.

\subsection{\label{subsec:comp_efficiency}Computational Efficiency}
In \cref{table:direct_second_order,table:direct_fourth_order,%
table:pw_second_order,table:pw_fourth_order} we give the timings for a
single time-step in the \textit{weak Landau damping} test-case. The hardware
hardware has the following specfications:
\begin{center}
\begin{tabular}{cc}
\toprule
CPU & Intel(R) Xeon(R) E-2276M  CPU @ 2.80GHz \\
    & 6 cores \\
\midrule
RAM & 32 GB DDR4 Synchronous 2667 MHz\\
\bottomrule
\end{tabular}
\end{center}

For the direct approach only resulutions up to $h_x =  L/64$, $h_v= \vmax/64$
were tested, due to memory constraints. When comparing the timings for a single
time-step for the direct and PW methods, we observe that the PW method is always
significantly faster. Even for low resolutions the PW method is still at least 
one order of magnitude faster than the direct method. Thus we conclude that from
a performance perspective the PW method  is indeed always the best choice. In
particular, \Cref{table:pw_second_order,table:pw_fourth_order} confirm that the
PW method scales linearly in the number of particles.

\Cref{fig:f_l2_error_2_order_kernel,fig:f_linfty_error_4_order_kernel} show that
for a given resolution, the direct method is up to two orders of magnitude more
accurate than the PW method. However, this comparison is misleading: the PW
requires higher resolutions to reach the same accuracy, but still needs
significantly less computational time to do so.

\begin{table}
\centering
\begin{tabular}{ccc}
\toprule
Resolution &  $t_{step}$ in s & $t_{total}$ in s \\
\midrule
$h_x = L/16$, $h_v= \vmax/16 $ & $2.03 \cdot 10^{-2}$ & 4.92 \\
$h_x = L/32$, $h_v= \vmax/32 $ & $1.78 \cdot 10^{-1}$ & 43.0 \\
$h_x =  L/64$, $h_v= \vmax/64$ & 5.62 & 1360 \\
\bottomrule
\end{tabular}
\caption{\label{table:direct_second_order}Timings for the direct approach 
using second order kernels.}
\end{table}

\begin{table}
\centering
\begin{tabular}{ccc}
\toprule
Resolution &  $t_{step}$ in s & $t_{total}$ in s \\
\midrule
$h_x = L/16 $, $h_v= \vmax/16 $ &  $2.43 \cdot 10^{-2}$ & 5.89 \\
$h_x =  L/32 $, $h_v= \vmax/32$ & $2.54\cdot 10^{-1}$ & 61.5 \\
$h_x = L/64$, $h_v= \vmax/64$ & 7.16 & 1730 \\
\bottomrule
\end{tabular}
\caption{\label{table:direct_fourth_order}Timings for the direct approach 
using fourth order kernels.}
\end{table}

\begin{table}
\centering
\begin{tabular}{ccc}
\toprule
Resolution &  $t_{step}$ in s & $t_{total}$ in s \\
\midrule
$h_x = L/16$, $h_v=\vmax/16$ & $1.22 \cdot 10^{-3}$  & 0.590 \\
$h_x = L/32$, $h_v= \vmax / 32$ & $2.77 \cdot 10^{-3}$ & 1.34 \\
$h_x = L/64$, $h_v=\vmax/64$ & $9.05 \cdot 10^{-3}$ & 4.28 \\
$h_x = L/128$, $h_v= \vmax / 256$ & $3.00 \cdot 10^{-2}$ & 14.3 \\
$h_x = L/256$, $h_v=\vmax/512$ & $1.15 \cdot 10^{-1}$ & 55.7 \\
\bottomrule
\end{tabular}
\caption{\label{table:pw_second_order}Timings for the PW approach using 
second order kernels.}
\end{table}

\begin{table}
\centering
\begin{tabular}{ccc}
\toprule
Resolution &  $t_{step}$ in s & $t_{total}$ in s \\
\midrule
$h_x = L/16$, $h_v=\vmax/16$ & $2.29 \cdot 10^{-3}$  & 1.11 \\
$h_x = L/32$, $h_v= \vmax / 32$ & $5.94 \cdot 10^{-3}$ & 2.88 \\
$h_x = L/64$, $h_v=\vmax/64$ & $1.90 \cdot 10^{-2}$ & 9.18 \\
$h_x = L/128$, $h_v= \vmax / 256$ & $6.45 \cdot 10^{-2}$ & 31.2 \\
$h_x = L/256$, $h_v=\vmax/512$ & $2.55 \cdot 10^{-1}$ & 123 \\
\bottomrule
\end{tabular}
\caption{\label{table:pw_fourth_order}Timings for the PW approach using 
fourth order kernels.}
\end{table}

\begin{table}
\centering
\begin{tabular}{ccc}
\toprule
Resolution &  $t_{step}$ in s & $t_{total}$ in s \\
\midrule
$h_x = L/16$, $h_v=\vmax/16$ & $1.07 \cdot 10^{-4}$  & $2.57 \cdot 10^{-2}$ \\
$h_x = L/32$, $h_v= \vmax / 32$ & $1.63 \cdot 10^{-4}$ & $3.92 \cdot 10^{-2}$ \\
$h_x = L/64$, $h_v=\vmax/64$ & $2.90 \cdot 10^{-4}$ & $6.94 \cdot 10^{-2}$ \\
$h_x = L/128$, $h_v= \vmax / 256$ & $9.37 \cdot 10^{-4}$ & $2.25 \cdot 10^{-1}$ \\
$h_x = L/256$, $h_v=\vmax/512$ & $4.40 \cdot 10^{-3}$ & 1.06 \\
\bottomrule
\end{tabular}
\caption{\label{table:PIC}Timings for a simple PIC method.}
\end{table}

\begin{figure}
\centering
\begin{subfigure}{0.48\textwidth}
\includegraphics[width=0.98\textwidth]{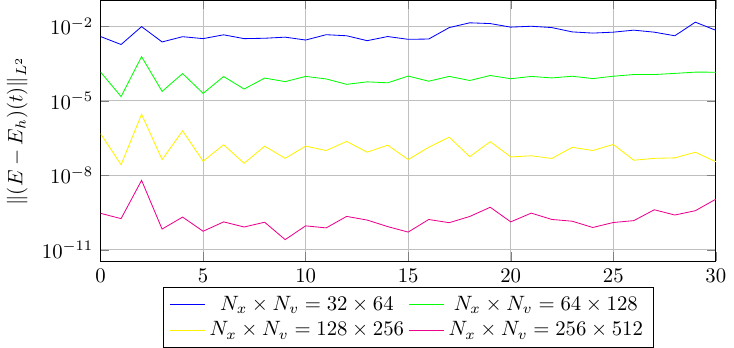}
\subcaption{\label{subfig:E_l2_error_4_order_kernel_pw} PW, 4.order}
\end{subfigure}
\begin{subfigure}{0.48\textwidth}
\includegraphics[width=0.98\textwidth]{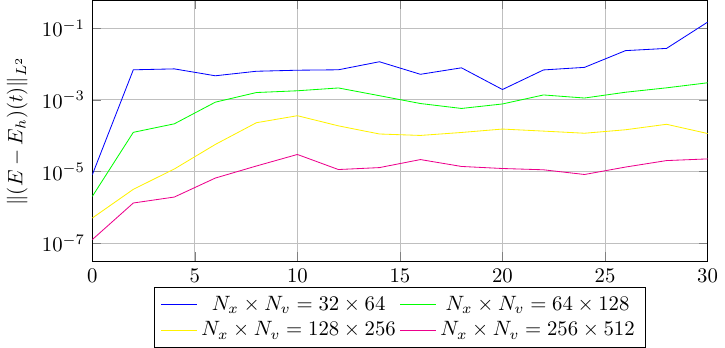}
\subcaption{\label{subfig:E_l2_error_PIC} PIC}
\end{subfigure}
\caption{\label{fig:E_l2_error_PIC_PW_4th_comparison} 
A comparison of $L^2$-errors for the numerically computed electric field
using on the one hand the PW method with fourth order kernels and 
$\sigma_x = 1$, $\sigma_v = 0.5$, and, on the other hand, a simple PIC 
method.}
\end{figure}

\begin{figure}
\centering
\begin{subfigure}{0.48\textwidth}
\includegraphics[width=0.98\textwidth]{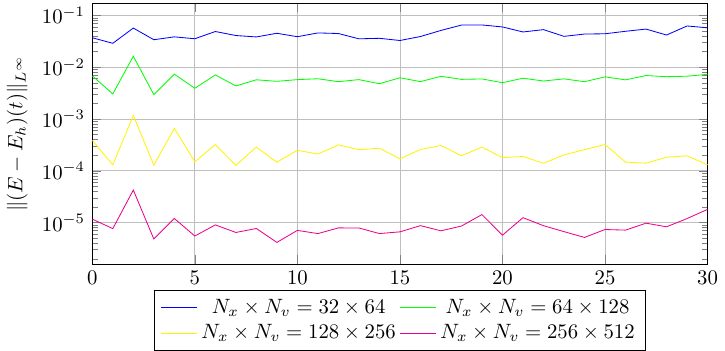}
\subcaption{\label{subfig:E_linfty_error_4_order_kernel_pw} PW, 4.order}
\end{subfigure}
\begin{subfigure}{0.48\textwidth}
\includegraphics[width=0.98\textwidth]{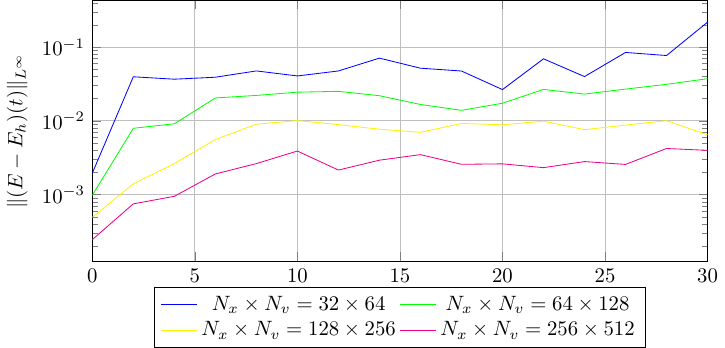}
\subcaption{\label{subfig:E_linfty_error_PIC} PIC}
\end{subfigure}
\caption{\label{fig:E_linfty_error_PIC_PW_4th_comparison} 
A comparison of $L^\infty$-errors for the numerically computed electric field
using on the one hand the PW method with fourth order kernels and 
$\sigma_x = 1$, $\sigma_v = 0.5$, and, on the other hand, a simple PIC 
method.}
\end{figure}

Finally we also compare the PW method with a simple PIC method. Our PIC
code approximates the density $\rho(x)$ for $jh_x\leq x<(j+1)h_x$ by adding
the masses of all particles in that $x$-range, and dividing the result
by $h_x$. The electric potential $\varphi$ is approximated using a standard
second-order finite-element method, time-integration uses the classical
Runge--Kutta scheme with time-step $\Delta t= \frac{1}{8}$.
\Cref{fig:E_l2_error_PIC_PW_4th_comparison,fig:E_linfty_error_PIC_PW_4th_comparison}
show the resulting errors in the elictric field.

One observes that the PW method converges significantly faster in both 
the $L^2$ and $L^\infty$-norm. This is in particular the case for the
finest resolutions. On the other hand, the simple PIC code is significantly
faster for a given resolution, as can be seen by comparing
\Cref{table:pw_fourth_order,table:PIC}. 
However, the high $L^\infty$-errors for the PIC method suggest
a strong level of numerical noise, which is much less present in the PW method.
Here, especially in the $L^2$-sense, it is much less clear at which level of
accuracy the PW method will begin to outperform the PIC method. On the other
hand, it is clear that the PIC method will require significantly larger numbers
of particles and thus imposes larger memory constraints on the machine.

\section{Conclusion}
We have presented a particle method using meshfree interpolation of arbitrary high order 
and investigated numerically, whether the good convergence behaviour from RKHS theory carries 
over to the case of the fully non-linear Vlasov-Poisson equation in the $d=1$ case. 

In contrast to conventional particle methods like PIC or SPH, our method does not need a 
remapping strategy to avoid numerical noise. Furthermore, as interpolation with Wendland
kernels is stable with arbitrary high convergence order, our method needs significantly
fewer particles to achieve the same accuracy as classical particle methods. 
The downside is that our method struggles with steep gradients, which naturally appear in
solutions of the Vlasov--Poisson equation. To resolve them correctly, any interpolation method 
needs high resolution, irrespective of the convergence order of the method. Similar problems
can be observed with higher order Eulerian Vlasov solvers. Steep gradients lead to overshoots.
However, while negative values of $f$ do not make any physical sense, their effect on the
quantities $\rho$ and $E$ seems to be limited.

The ill-conditioning of the kernel matrices is a well-known problem and its solution is an
ongoing research topic in the RKHS community. This limits the particle numbers for the direct 
method. In case of the Vlasov--Poisson equation, however, this problem can to some extend
be bypassed by using piece-wise interpolants instead. Our simulations have illustrated that 
this approach does in fact result in efficient Lagrangian schemes, albeit with convergence
orders that are lower than those of the direct approach.

To summarise, at least in the one-dimensional case, the presented PW method offers a
good compromise between the stability and  high accuracy of purely Eulerian methods on the
one hand and speed and hyperbolicity of classical particle methods on the other hand, while
avoiding the inherent numerical noise of the latter. 

It is unclear whether the piecewise approach is suitable for higher dimensions, for which
we expect reduced efficiency as larger local systems need to be solved. On the other hand,
we believe that the piece-wise approach could prove to be beneficial for stellar dynamics,
where particles tend cluster more stronlgy, and could thus provide \enquote{auto-adaptation}.
These points require further investigation.

\printbibliography
	
\end{document}